%% file: main.tex
\documentclass[letterpaper, 10 pt, conference, dvipsnames]{ieeeconf}  

\IEEEoverridecommandlockouts                              

\usepackage{mystyle}

\usepackage{notation}
\usepackage{scenariotree}
\usepackage[nameinlink,capitalise]{cleveref}
\usepackage{amsmath} 
\usepackage{amssymb}  
\usepackage{bbm}

\title{Interaction-aware Model Predictive Control for Autonomous Driving
}

\author{Renzi Wang \and Mathijs Schuurmans \and Panagiotis Patrinos%
\thanks{This work was supported by
the Research Foundation Flanders (FWO) research projects G081222N, G033822N, G0A0920N;
European Union’s Horizon 2020 research and innovation programme under the Marie Skłodowska-Curie grant agreement No. 953348;
Ford KU Leuven Research Alliance Project KUL0075; 
Research Council KU Leuven C1 project No. C14/18/068;
Fonds de la Recherche Scientifique-FNRS the Fonds Wetenschappelijk Onderzoek-Vlaanderen under EOS project no 30468160 (SeLMA).}%
\thanks{
    KU Leuven, Department of Electrical Engineering \textsc{esat-stadius} -- %
			Kasteelpark Arenberg 10, bus 2446, B-3001 Leuven, Belgium
    \newline
        {\sf
            \{%
                \href{mailto:renzi.wang@kuleuven.be}{renzi.wang},
                \href{mailto:mathijs.schuurmans@kuleuven.be}{mathijs.schuurmans},
                \href{mailto:panos.patrinos@kuleuven.be}{panos.patrinos}%
            \}%
            \href{mailto:renzi.wang@kuleuven.be,mathijs.schuurmans@kuleuven.be,panos.patrinos@kuleuven.be}{@kuleuven.be}%
    }
}
}

\begin{document}

\maketitle
\thispagestyle{empty}
\pagestyle{empty}

\input{content/abstract}


\section{INTRODUCTION}
\input{content/introduction}



\section{Modeling and Problem Formulation}\label{sec: formulation}

\input{content/problem_formulation}

\section{Interaction-Aware MPC}\label{sec: ineraction-mpc}

\input{content/interaction_aware_mpc}



\section{Numerical experiments}\label{sec: experiment}

\input{content/experimental_results}

\section{Conclusions}

\input{content/conclusion_and_future_work}\label{sec: conclusion}









\bibliographystyle{IEEEtran}
\bibliography{IEEEabrv,reference}

\end{document}

%% file: content/abstract.tex
\begin{abstract}
    Lane changing and lane merging remains a challenging task for autonomous driving, due to the strong interaction between the controlled vehicle and the uncertain behavior of the surrounding traffic participants. 
    The interaction induces a dependence of the vehicles’ states on the (stochastic) dynamics of the surrounding vehicles, increasing the difficulty of predicting future trajectories. 
    Furthermore, the small relative distances cause traditional robust approaches to become overly conservative, necessitating control methods that are explicitly aware of inter-vehicle interaction. 
    Towards these goals, we propose an interaction-aware stochastic model predictive control (MPC) strategy integrated with an online learning framework, which models a given driver's cooperation level as an unknown parameter in a state-dependent probability distribution. 
    The online learning framework adaptively estimates the surrounding vehicle’s cooperation level with the vehicle’s past trajectory and combines this with a kinematic vehicle model to predict the probability of a multimodal future state trajectory. 
    The learning is conducted with logistic regression which enables fast online computation. 
    The multi-future prediction is used in the MPC algorithm to compute the optimal control input while satisfying safety constraints. 
    We demonstrate our algorithm in an interactive lane changing scenario with drivers in different randomly selected cooperation levels.
\end{abstract}

%% file: content/introduction.tex
The development of autonomous vehicles has been a long-standing endeavour of 
research and industry due to its potentially far-reaching implications 
on safety and efficiency of our traffic.
Despite significant advances towards this goal, one of the key remaining challenges
is the problem of control and planning in close 
vicinity to surrounding traffic participants.
The main sources for these challenges are
\begin{inlinelist*}
    \item the inherent uncertainty about the intentions of surrounding drivers
    \item the interaction between the behavior of these drivers and that of the controlled vehicle. 
\end{inlinelist*}
To address these issues, we propose a stochastic model for human driving behavior in 
a highway driving scenario, 
in which the human driver is assumed to repeatedly make a selection 
out of a finite set of base policies according to a 
(\textit{a priori} unknown) probability distribution. In order to model interactive behavior, we 
will pose the natural assumption that this distribution depends on the joint state of the involved vehicles 
(e.g., relative positions and velocities).
We describe a simple, online scheme for estimating the dependence 
of the aforementioned probability distribution on the 
state and include this model in a chance-constrained stochastic \ac{MPC} formulation, which 
can be solved effectively over a scenario tree. 

\subsection{Related Work}
\subsubsection{Modeling uncertainty}
Existing approaches to \ac{MPC} under uncertainty on surrounding vehicle motion
have can be broadly categorized into \textit{robust} \cite{bemporad1999robust}
and \textit{stochastic} approaches,
where the former assumes potentially adversarial behavior of surrounding vehicles, which 
is accounted for under \textit{worst-case} assumptions and the latter 
assumes a stochastic description of the driving behavior, where 
the expected value of the cost is minimized subject 
to potentially probabilistic constraints \cite{mesbah2016stochastic}.
More recently,
risk-averse \cite{sopasakis2019risk,sopasakis2019risk2,schuurmans2020learningacc,dixit_RiskAverseRecedingHorizon_2022} and distributionally robust \cite{schuurmans2019safe,schuurmans2020learningcdc,coppens2020data,schuurmans2021data,coppens2021data,schuurmans_SafeLearningBasedMPC_2022b,schuurmans2021general}
approaches have become increasingly popular as they generalize the aforementioned approaches,
by accounting for uncertainty in the distribution of the stochastic prediction model.

Most of these existing methods, however, consider the uncertainty to be fully exogenous --- the probability distribution
does not depend on the state or inputs of the controlled vehicle. 
This tends to result in reactive and overly conservative behavior (cf. the \textit{frozen robot problem} \cite{trautman_UnfreezingRobotNavigation_2010}). 

\subsubsection{Modeling interaction}
A variety of approaches has been proposed to include the interactive behavior in the controller design.
Reinforcement learning is applied in \cite{bouton2020reinforcement,hu2019interaction} to learn an interaction-aware policy directly from data. 
These methods require little expert knowledge and can typically handle complicated scenarios. However, 
the resulting policies are usually not easy to introspect, making it difficult to ensure safety or predict 
failure cases. Furthermore, their quality hinges on the availability 
of large-scale high-quality datasets.

An alternative approach is to combine \ac{MPC} with predictive models that describe inter-vehicle interaction.
The benefit of this approach is that \ac{MPC} simultaneously allows for a wide range of expressive dynamical models 
and furthermore directly accounts for safety restrictions through constraints in the optimal control problem.

Several different interaction-aware \ac{MPC} formulations have been proposed in recent years, 
such as game-theoretical \ac{MPC} \cite{dreves_GeneralizedNashEquilibrium_2018a,fabiani2019multi,cleac2019algames,evens2022learning}, 
where each road user is assumed to behave optimally according to some known or estimated cost function. 

A more direct approach is to learn the interactive dynamics through supervised learning \cite{bae2020cooperation, liu2021interaction}, 
where given current states of each vehicle, a trained neural network predicts their future trajectories in a prediction time horizon.
Both methods considers uni-modal prediction where the prediction is based on most likely scenario or a predefined policy. 

As mentioned before, our goal is to additionally account for the inherent uncertainty in the driving behavior, by considering a stochastic interaction 
model. Our model class is similar to \cite{chen_InteractiveMultiModalMotion_2022}. However, in \cite{chen_InteractiveMultiModalMotion_2022}, 
the (state-dependent) distribution is assumed to be known, whereas we integrate a procedure for 
adapting a model for these distributions during operation.
Simiarly, \cite{hu2022active} integrates \ac{MPC} with an online learning technique for a similar use case.
However, in \cite{hu2022active}, 
human driving behavior is modeled as a weighted sum (with unknown weights) of random policies with known distributions, which depends on an estimate of the human driver's action-value function (i.e., their long-term costs as a function of the current state and control action).
Instead, we assume that the human driver will select one of the predefined number of base policies according to an unknown, state-dependent distribution. This only requires knowledge of
these base policies corresponding to possible maneuvers, which provides a very direct and intuitive way of modeling possible behaviors for a given traffic scenario.
Furthermore, our formulation allows to reduce excessive conservatism through chance constraints, which 
allows to discard extremely unlikely scenarios.  
However, \cite{hu2022active} presents an indirect dual adaptive control method for active learning, which could be 
applicable to the model considered in our work. We leave a more thorough investigation of this possibility
for future work. 

\subsection{Contribution}
We tackle a lane change control problem by utilizing \ac{SMPC} \cite{patrinos2014stochastic}, 
where we optimize the control input under a finite number of scenarios under decision-dependent uncertainty. 
Each scenario is defined with a maneuver-based prediction model \cite{lefevre2014survey}.
In the machine learning literature, this model class is also known as \textit{multi-future prediction} \cite{chai2019multipath, cui2019multimodal}.
Different from \cite{cesari2017scenario, nair2022stochastic}, we explicitly consider interaction by introducing a state-dependent distribution.
Our contribution is twofold:
\begin{inlinelist*}
    \item We integrate our \ac{MPC} method with a simple learning framework which, during operation, adapts the probabilistic interaction model to the current neighboring driver.
    \item We present \iac{SMPC} method with a smooth, safe approximation of the (originally non-smooth) chance constraints.
\end{inlinelist*}


\subsection*{Notation}
Let $\norm{x}_2$ denote the Euclidean norm of vector $x$ and $\norm{x}_A^2 \dfn x^\top\!Ax$ for a positive definite matrix $A\succ 0$.
We use $\diag{x}$ to represent a diagonal matrix with diagonal elements equal to the entries of vector $x$. 
Given a set $C \subseteq \Re^n$, $\indi[C]: \Re \ni x \mapsto \indi[C](x) \in \Re$ denotes the function that maps $x$ to 1 when $x \in C$ and to 0 otherwise.
$U[a, b]$ represents the uniform distribution defined on the closed interval $[a, b] \subset \Re$ .
$\mathcal{N}(\mu, \Sigma)$ denotes the normal distribution with mean $\mu$ and covariance matrix $\Sigma$.
Furthermore, for $a < b \in \N$, we denote $\natseq{a}{b} \dfn \{k \mid a \leq k \leq b\}$.
We denote by $\simplex_{\nModes} \dfn \{ p \in \Re^{\nModes}_+ \mid \sum_{i=1}^{\nModes} p_i = 1 \}$ the 
$\nModes$-dimensional probability simplex.
We denote by $\|A\|_\fro$ the Frobenius norm of the matrix $A$.

%% file: content/problem_formulation.tex
Consider a traffic scenario on a highway, consisting of a controlled 
vehicle (the \textit{ego vehicle}), referred to by the identifier $\ego$, 
and a number of surrounding road vehicles.
The task of the ego vehicle is to safely perform a lane-change 
maneuver to a lane which is occupied by other road users. 
At the time of initiating the lane change procedure, 
we consider the nearest vehicle approaching from the rear 
in the adjacent lane, which we refer to as the 
\textit{target vehicle}, and identify with 
the index $\tv$.
Since the behavior of this agent determines 
whether the lane change can be safely carried out, 
we will explicitly model the interaction dynamics 
with this vehicle, and assume that remaining obstacles 
can be handled using simplified models. 

\subsection{Vehicle Dynamics}
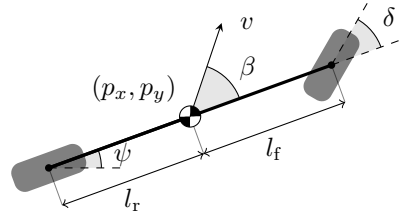
\begin{figure}[tb]
    \centering
    \input{tikz/bicycle_drawing.tex}
    \caption{Illustration of the kinematic bicycle model.}
    \label{fig:kinematic-bicycle-model}
\end{figure}

We model all vehicle dynamics with a kinematic bicycle model~\cite{polack2017kinematic}, 
illustrated in \cref{fig:kinematic-bicycle-model}. 
Let $\pos = \smallmat{\pos_x & \pos_y}^\top \in \Re^2$ describe the position of the center of gravity expressed in a fixed inertial frame,
$\psi$ is the orientation of the vehicle, and $v$ is the velocity of the center of gravity.
Its control inputs are the acceleration $a$ and the steering angle $\steer$. 
The continuous time dynamics are given as 
\begin{subequations}\label{eq: kinematic_bicycle}
\begin{align}
    \dot{\pos}_x =&\; v\cos(\psi + \beta) \\
    \dot{\pos}_y =&\; v\sin(\psi + \beta) \\
    \dot{v} =&\; a \\
    \dot{\psi} =&\; \tfrac{v}{\lr} \sin \beta \\
\intertext{where $\beta$ is the slip angle at the center of gravity:}
    \beta =&\; \arctan\Big(\tfrac{\lr}{\lf + \lr}\tan \steer \Big).
\end{align}
\end{subequations}
Denoting the state of vehicle $z \dfn \smallmat{\pos_x & \pos_y & v & \psi}^\top \in \Re^4$ and 
the control action $u \dfn \smallmat{a & \steer}^\top \in \Re^{\na}$ ($\na = 2$), 
we discretize \eqref{eq: kinematic_bicycle} using the forward Euler method
to obtain a discrete-time system consistent with the simulator we used in the experiments, 
which we compactly denote as
\begin{equation}\label{eq:dt-dynamics}
    z_{k+1} = f(z_k, u_k).
\end{equation}
For ease of notation, we will denote the state of the ego vehicle and the target vehicle as 
\(\idx{z}{\ego}\) and \(\idx{z}{\tv}\), respectively, and we define the 
concatenated state $\full{z} \dfn \left( \idx{z}{\ego}, \idx{z}{\tv} \right) \in \Re^{\ns}$ (with $\ns = 8$).

\subsection{Target vehicle behavior} \label{sec: finite_scenarios}
Let $(\Omega, \mathcal{F}, \prob)$ be a probability space with 
countably infinite sample space $\Omega$, $\sigma$-algebra $\mathcal{F}$
and probability measure $\prob: \mathcal{F} \to [0,1]$.
We model the uncertain behavior of a given target vehicle
as a randomized control policy $\idx{\kappa}{\tv}: \Re^{\ns} \times \W \to \Re^2$, 
where $\W \dfn \{1, \dots, \nModes\}$ is a finite set corresponding to 
a number of discrete maneuvers that a driver may perform.
We assume that every time step $t$,
the selection of the maneuver $\md_t \in \W$ 
    is made randomly according to a distribution 
\begin{equation} \label{eq:def-ddprob}
    \ddP\left( \full{z}_t \right) = (\prob\{ \md_{t} = i \mid \full{z}_t \})_{i \in \W}, 
\end{equation}

which depends on the (combined) state $\full{z}_t$ of the vehicles at time $t$.
This state-dependent distribution $\ddP: \Re^{\ns} \to \simplex_\nModes$ will typically differ  
between different drivers, and can therefore not easily be estimated 
offline using traditional methods. 
Instead, we propose a simple online estimation scheme to estimate 
this distribution online; This is described in \cref{sec: learning}.
Thus, at a given time $t$, the control action of the target vehicle is a 
random variable
\[ 
    \idx{u}{\tv}_{t} = \kappa^\tv \left( \full{z}_t, \md_{t} \right), \quad \md_{t} \sim \ddP(\full{z}_t), 
\]
which, combined with \eqref{eq:dt-dynamics} results in an autonomous, stochastic system
\(
    \idx{z}{\tv}_{t+1} = \idx{f}{\tv}\left( \full{z}_t, \md_t \right) \dfn f\left( \idx{z}{\tv}_t , \idx{\kappa}{\tv}\left( \full{z}_t, \md_t \right)   \right), 
\)
describing the target vehicle motion.
For ease of notation, we combine the 
dynamics of the ego vehicle and the target vehicle, and introduce $\fulldyn: \Re^{\ns} \times \Re^{\na} \to \Re^{\ns}$
\begin{equation} \label{eq:full-dynamics}
    \full{z}_{t+1} = \fulldyn(\full{z}_t, u_t, \md_t) 
        \dfn \smallmat{ 
            f\left( \idx{z}{\ego}_t, u_t \right) \\ 
            \idx{f}{\tv}(\full{z}_t, \md_t)
        }.
\end{equation}
\begin{inlinelist*}
\item The model uncertainty is concentrated in the behavior, resulting in 
a restriction to dynamically feasible predicted trajectories
\item The choice of considering a fixed number of known policies 
allows the designer to incorporate a large amount of known 
structure about the use case at hand into the target vehicle's 
behavioral model
\item The model is similar to the popular technique of using \textit{motion primitives}\cite{lefevre2014survey}, but by using a control policy for the closed-loop model, 
the predicted trajectories automatically get corrected online using feedback
\item In contrast to similar models in which the distribution of $(\md_t)_t$ is state-independent, e.g.,
\cite{schuurmans_SafeLearningBasedMPC_2022b}, the proposed model grants considerably more 
freedom in modeling interaction between vehicles, since the probability of choosing a certain 
maneuver typically depends on the relative positions and velocities of nearby vehicles.
\end{inlinelist*}

\subsection{Objective Function}\label{sec: objective}
We use the usual quadratic stage cost $\ell: \Re^{\ns} \times \Re^{\na} \to \Re_+$,
\begin{equation}\label{eq:stage-cost}
    \ell(\full{z}, u) \dfn \|\idx{z}{\ego} - z^{\refer}\|^2_Q + \|u^\ego\|^2_R,
\end{equation}
where $Q \succeq 0$ and $R \succ 0$ are given weights for the ego vehicle's states and inputs, respectively, 
and $z^{\refer}$ is a given reference state. 

\subsection{Input and state constraints}
\subsubsection{Simple state and input bounds}
As usual in control applications, we will assume there are hard limits 
on both the states and the inputs, which we will 
write as 
\[
    u_{t} \in \U, \full{z}_t \in \Z, \forall t \in \N,
\]
where $\U \subseteq \Re^{\na}, \Z \in \Re^{\ns}$ are given, 
nonempty, closed sets. Since we are only actively controlling the 
ego vehicle, we will typically take $\Z = \idx{\Z}{\ego} \times \Re^4$, 
with $\idx{\Z}{\ego}$ a given constraint set for the ego vehicle state. 
In practice, both $\idx{\Z}{\ego}$ and $\U$ will often be simple
component-wise upper and lower bounds.

\subsubsection{Slew rate constraints}
To ensure comfortable driving behavior and account for unmodeled actuator dynamics, 
it is common to include slew-rate constraints of the form 
\[ 
 | u_{k} - u_{k-1} | \leq \dumax, 
\]
where the absolute values and inequality are taken element-wise, and $\dumax > 0$
is a given constant. 
\subsection{Collision Avoidance}
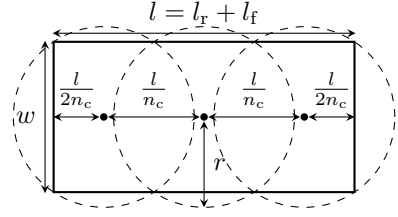
\begin{figure}[tb]
    \centering
    \newcommand\ncirc{2}
    \input{tikz/collision_avoidance}
    \caption{Illustration of the collision avoidance constraints.}
    \label{fig: collision_avoidance}
\end{figure}

Given a vehicle with current state $z = \smallmat{\pos & \vel & \head}^\top$,
with $\pos \in \Re^{2}$,
we -- similarly to, e.g., \cite{wang_NonGaussianChanceConstrainedTrajectory_2020} --
approximate the area it occupies as the union of 
$\ncircle = 3$ circles with radius $r = \tfrac{1}{2}\sqrt{(\nicefrac{l}{\ncircle})^2 + w^2}$ and centers
\begin{equation*} 
    c_{j}\left( z \right) \dfn \pos + \tfrac{l}{2\ncircle} \left( 2j - \ncircle -1 \right) \smallmat{\cos \head \\ \sin \head}  , \quad j = 1, \dots, \ncircle,
\end{equation*}
as illustrated in \cref{fig: collision_avoidance}.

A sufficient condition for collision avoidance between 
the vehicles $\ego$ and $\tv$, with respective states $z^\ego$ and $z^\tv$,
is then
\begin{equation}\label{eq: collision_avoidance_robust}
	g_{ij}\left( \full{z} \right) \dfn 4 r^2 - \| c_{i}(z^{\ego}) - c_{j}(z^{\tv}) \|_2^2 \leq 0, \quad \forall i,j \in [\ncircle].
\end{equation}
Note that the functions $g_{ij}: \Re^{\ns} \to \Re$ smooth, yet concave, and will 
thus lead to nonconvex constraints.

Very often, the distribution $\ddP\left( \full{z} \right)$ will be concentrated 
on a small subset of $\W$. For instance, whenever the ego vehicle drives 
directly in front of the target vehicle. The probability that the target vehicle 
will heavily accelerate, causing a rear-end collision, will be very low. 
Imposing collision avoidance constraint \eqref{eq: collision_avoidance_robust}
for the predicted states under this scenario, would be unreasonably conservative 
and would in many cases simply lead to infeasibility. 
To resolve this issue, we will instead aim to impose \eqref{eq: collision_avoidance_robust}
probabilistically. More specifically, at any time step $t$, our goal is to impose 
\[
    \prob\{ \max_{i,j \in [\ncircle]} g_{ij}\left( \full{z}_{t+1} \right) \leq 0 \mid \full{z}_t \} \geq 1 - \gamma, 
\]
for a given violation rate $\gamma \in [ 0, 1 )$. By 
definition \eqref{eq:def-ddprob}, this constraint depends on the state-dependent maneuver distribution $\ddP$.


\subsection{Nominal optimal control problem}
We are now ready to formulate the stochastic optimal control 
problem, which we ideally would like to solve 
in a receding horizon fashion to obtain an interaction-aware 
\ac{MPC} scheme for the described task. As we will see, however,
this problem cannot be solved directly. Instead, \cref{sec: ineraction-mpc} 
is dedicated to proposing suitable approximations to the formulations, resulting 
in a practical control scheme. 

Let 
\(
\Pol \dfn \{ \pol_k: \Re^{\ns} \to \Re^{\na}, k \in \natseq{0}{N-1} \}
\)
denote the set of all causal $N$-step control policies.

Given an initial state 
\(
\overline{\full{z}}
\), combine the introduced elements to obtain 
obtain the (idealized) stochastic optimal control problem
\begin{subequations}\label{eq: formulation_origin}
\begin{align}
	&\minimize_{(\pol_t)_{t} \in \Pol} && \E_\md \Big\{ \sum_{k = 0}^{N-1} \ell(\full{z}_k, u_k) + \lN\left( \full{z}_N \right)  \mid \bar{\full{z}} \Big\} \label{eq: obj_origin}\\
	&\;\stt&& \full{z}_0 = (\idx{z}{\ego}_0, \idx{z}{\tv}_0) = \bar{\full{z}},\\
	&&&z^\ego_{k+1} = f(\idx{z}{\ego}_k, u_k), u_{k} = \mu_{k}\left( \full{z}_k \right),
        \\
	&&&z^\tv_{k+1} = \idx{f}{\tv}\big(z^\tv_k, \md_{k}\big),
        \\
    &&&\prob
        \big\{ 
            \max_{i,j \in [\ncircle]} g_{ij}(\full{z}_{k+1}) 
                {}\leq{}
            0 \mid \full{z}_k 
        \big\} 
            \geq 1 - \gamma, \label{eq:ocp-chance-constraint}
        \\
    &&&u_k \in \U, \; \full{z}_k \in \Z, \full{z}_N \in \Z, \, \mathrm{a.s.},
        \\
    &&& | u_{k'+1} - u_{k'} | \leq \dumax \, \mathrm{a.s.}, \\ \notag
    &&&\forall k \in \natseq{0}{N-1},\; \forall k' \in \natseq{0}{N-2},
\end{align}
\end{subequations}
where the qualifier $a.s.$ (almost surely) indicates that the constraint must hold 
in all realizations that occur with nonzero probability.
Of course, as it is stated here, this problem is intractable due to 
\begin{inlinelist*}
    \item the infinite-dimensional decision variables $(\pol_k)_k$
    \item the distribution $\ddP$ governing the stochastic process $(\md_{k})_k$ (cf. \eqref{eq:def-ddprob}) being unknown 
    \item the non-smooth chance constraint \eqref{eq:ocp-chance-constraint}.
\end{inlinelist*}
In the next section, we will propose an approximate formulation of \eqref{eq: formulation_origin}, addressing these difficulties 
and resulting in an implementable control scheme.

%% file: tikz/bicycle_drawing.tex
\usetikzlibrary{calc}


    \newcommand\centerofmass{%
        \tikz[radius=0.4em] {%
            \fill (0,0) -- ++(0.4em,0) arc [start angle=0,end angle=90] -- ++(0,-0.8em) arc [start angle=270, end angle=180];%
            \fill [color=white] (0,0) -- ++(0,0.4em) arc [start angle=90,end angle=180] -- ++(0.8em,0) arc [start angle=0, end angle=-90];%
            \draw (0,0) circle;%
        }%
    }
    \begin{tikzpicture}
        \def\headangle{20}
        \def\steerangle{40}
        \def\tirelen{1}
        \def\tirewidth{0.2}
        \def\lfdist{2}
        \def\lrdist{2}
        \def\ticklen{15pt}

        \tikzset{tire/.style={rounded corners, fill=gray}, 
                 car/.style={very thick}, 
                 auxline/.style={thin, dashed},
                 point/.style={circle, inner sep=1, fill=black},
                 angle/.style={fill=black, fill opacity=0.1, draw opacity=1, draw=black},
                 tick/.style={very thin, draw opacity=0.5}, 
                 arrow/.style={<->, thin}
                }

        \coordinate (back) (1,1);
        \fill[rotate={\headangle}, tire] ($(back) - (0.5,0.2)$) rectangle ($(back) + (0.5,0.2)$);
        \node[point] at (back) {};
        \draw[rotate={\headangle}, car, anchor=west] (back) -- ($(back) + ({\lrdist}, 0)$) coordinate (com);
        \path[rotate={\headangle}] (back) -- ($(com) + ({\lfdist}, 0)$) coordinate (front);
        \fill[rotate={\headangle+\steerangle}, tire] ($(front) - (0.5,0.2)$) rectangle ($(front) + (0.5,0.2)$);
        \draw[rotate={\headangle}, car, anchor=west] (back) -- ($(com) + ({\lfdist}, 0)$);

        \draw[rotate={\headangle}, anchor=west, tick] (com) -- ++(0, -15pt) coordinate (nodemiddle);
        \draw[rotate={\headangle}, anchor=west, tick] (back) -- ++(0, -15pt) coordinate (nodeback);
        \draw[rotate={\headangle}, anchor=west, tick] (front) -- ++(0, -15pt) coordinate (nodefront);
        
        \draw[rotate={\headangle}, car, anchor=(back), arrow] (nodeback) -- (nodemiddle) node [midway, below] {$\lr$};
        \draw[rotate={\headangle}, car, anchor=(back), arrow] (nodemiddle) -- (nodefront) node [midway, below] {$\lf$};
        \node[label={[label distance=-0.3cm]above left:{$(\pos_{x}, \pos_{y})$}}] at (com) {\centerofmass};
        \pgfmathparse{atan(0.5*tan(\steerangle/180*pi))}
        \draw[->,>=stealth, rotate={\headangle + \pgfmathresult*180/pi}] (com) -- (3, 0) node[label=right:{$v$}](velocity){};
        \node[point] at (front) {};
        \begin{scope}
            \path[clip] (front) -- (com) -- (velocity);
            \draw[angle] (com) circle (0.7);
        \end{scope}
        \node at ($(com) + ({\headangle + (\pgfmathresult*180/pi)}:1)$) {$\beta$};

        \draw[auxline] (back) -+ (1, 0) coordinate (auxback);
        \begin{scope}
            \path[clip] (auxback) -- (back) -- (front);
            \draw[angle] (back) circle (0.7);
        \end{scope}
        \node at ($(back) + ({0.5*\headangle}:1)$) {$\psi$};

        \draw[auxline, rotate={\headangle}, anchor=west] (front) -- ($(front) + (1, 0)$) coordinate (auxfront);
        \draw[auxline, rotate={\headangle+\steerangle}, anchor=west] (front) -- ($(front) + (1, 0)$) coordinate (auxfrontwheel);
        \begin{scope}
            \path[clip] (auxfrontwheel) -- (front) -- (auxfront);
            \draw[angle] (front) circle (0.7);
        \end{scope}
        \node at ($(front) + ({\headangle + 0.5*\steerangle}:1)$) {$\delta$};

    \end{tikzpicture}

%% file: tikz/collision_avoidance.tex
\begin{tikzpicture}
    \tikzset{point/.style={circle, inner sep=1pt, fill}}
    \providecommand{\distsym}{\tfrac{l}{2\ncircle}}
    \providecommand{\distsymdouble}{\tfrac{l}{\ncircle}}
    \newlength{\w}\setlength{\w}{4cm}
    \newlength{\he}\setlength{\he}{2cm}
    \node[draw=black, thick, minimum width=\w, minimum height=\he] (rect){};
    \draw[<->, >=stealth] ($(rect.north west) + (0,0.1cm)$) -- ($(rect.north east) + (0,0.1cm)$) node[midway, above] {$l = \lr + \lf$};
    \draw[<->, >=stealth] ($(rect.north west) - (0.1cm,0)$) -- ($(rect.south west) - (0.1,0cm)$) node[midway, left] {$w$};
    \providecommand\ncirc{1}
    \pgfmathsetmacro{\dl}{1/(2*(\ncirc+1))*\w}
    \pgfmathsetmacro{\circrad}{\he / (2 * sin(atan(\he*(\ncirc+1)/\w)))}
    \coordinate (c--1) at (-0.5\w, 0);
    \foreach\i in {0, 1, ..., \ncirc}{
        \pgfmathparse{(1 + 2*\i)*\dl};
        \node[point] (c-\i) at (-0.5\w + \pgfmathresult pt, 0) {};
        \draw[dashed] (c-\i) circle (\circrad pt);
    }
    \draw[<->, >=stealth] (c--1) -- (c-0) node[midway, above] {$\distsym$};
    \foreach\i in {1, ..., \ncirc}{
        \pgfmathparse{\i-1}
        \draw[<->, >=stealth] (c-\pgfmathresult) -- (c-\i) node[midway, above] {$\distsymdouble$};
    }
    \draw[<->, >=stealth] (c-\ncirc) -- (0.5\w, 0) node[midway, above] {$\distsym$};
    \pgfmathparse{}
    \draw[<->, >=stealth] (c-1) -- ++(0, -\circrad pt) node[midway, right] {$r$};
\end{tikzpicture}

%% file: content/interaction_aware_mpc.tex
\subsection{Learning the distribution}\label{sec: learning}
\input{content/learning_the_distribution}

\subsection{Chance constraint approximation} \label{sec: chance-constraints}
\input{content/tractable_reformulations}

%% file: content/learning_the_distribution.tex
The state-dependent distribution $\ddP$ introduced in \eqref{eq:def-ddprob} 
is not known in advance, and furthermore, 
it cannot be accurately estimated offline,
as it may vary significantly between drivers.
This necessitates the introduction of an \textit{online} learning method. 

To this end, we use a simple 
multinomial logistic regression scheme,
where we introduce a parameterized model $\hat{\ddP}: \Re^{\ns} \times \Re^{\nfeat \times \nModes} \to \Re^{\nModes}$
\begin{equation} \label{eq:model-probability}
    \hat{\ddP}\left( \full{z} ; \param  \right) \dfn
    \bigg( 
        \frac{\exp(\param_i^\top \feat(\full{z}))}{\sum_{j=1}^{\nModes}\exp(\param_i^\top \feat(\full{z}))}
    \bigg)_{i = 1}^{\nModes}
\end{equation}
where $\feat: \Re^{\ns} \to \Re^{\nfeat}: \full{z} \mapsto \left( \feat_{i}\left( \full{z} \right)  \right)_{i = 1}^{\nfeat}$,
 is a vector-valued feature function, with $\feat_1 \equiv 1$, which accounts for the bias term, and 
$\param = \smallmat{\param_{1} & \dots & \param_{\nModes}} \in \Re^{\nfeat \times \nModes}$ is the 
parameter matrix.


\subsubsection{Maximum-likelihood estimation}
Given a dataset $D_T \dfn \left( \full{z}_k, \md_{k}  \right)_{k = 1}^{T}$, 
the parameter $\param$ is determined by maximizing the log-likelihood \cite[\S 4.4]{hastie2009elements}:
\begin{equation}\label{eq:offline-param}
    \param^\star \in \argmin_\param -\tsum_{k=1}^{T} \log \hat{\ddP}_{\md_k}\big(\feat(\full{z}_k);\param\big).
\end{equation}
If $D_T$ is obtained offline from a collection of different drivers, then it may be used 
to obtain an initial guess for the driving parameters of a new driver. However, 
since $\param$ may vary significantly from one driver to another, 
we additionally consider an online procedure for updating $\param$ for the 
currently neighboring target vehicle.

\subsubsection{Online Learning}\label{sec: online_learning}
In order to adapt the parameter with measurements obtained 
during operation, we use a moving horizon estimation update. 
At every time step $t$, we compute a new estimate
\begin{equation}\label{eq: param_update}
    \param_{t} \in \argmin_{\param} \lambda \reg(\param, \param_{t-1}) - \sum_{k=t-L+1}^t \log \hat{\ddP}_{\md_k} \big(\full{z}_k; \param \big),
\end{equation}
where $L \in \N$ is a fixed window length, $\lambda > 0$ is a fixed regularization constant and 
$\reg: \left( \param, \param' \right) \mapsto \| \param - \param' \|_{\fro}^2$
penalizes large differences between $\param$ and the previous estimate $\param_{t-1}$. 


%% file: content/tractable_reformulations.tex
We now move our attention to obtaining approximate reformulations of the 
problem \eqref{eq: formulation_origin}, which can be 
effectively solved online within \iac{MPC} scheme.

Let us first consider the chance constraint  
\eqref{eq:ocp-chance-constraint}.
Given a random quantity $\zeta: \W \to \Re$, we may write 
\begin{equation} \label{eq:chance-exact}
    \prob\{ \zeta(\md_k) > 0 \mid \full{z}_k \} = \sum_{\md \in \W} \ddP_{\md}(\full{z}_k) \indi[{\Re_+}](\zeta(\md_k)).
\end{equation}
The function $\indi[{\Re_+}]$ is nonconvex, and furthermore, it is discontinuous at 0, making it unsuitable 
for numerical optimization. To remedy this, we aim to replace it by a continuous and/or smooth 
approximation. A common choice for such a surrogate function is the average value-at-risk ($\AVAR$) \cite{nemirovski_SafeTractableApproximations_2012}, 
which results in the smallest convex upper approximation of the original chance constraint. However, it may be rather conservative, 
in particular in the case where the number of outcomes, i.e., $|\W| = \nModes$, is small, as we illustrate in \cref{ex:chance-constraints}. 
Due to the collision avoidance constraints and the nonlinear dynamics, \eqref{eq: obj_origin} 
is nonconvex
irrespective of 
the surrogate used for $\indi[{\Re_+}]$.
Obtaining a convex upper approximation is therefore of limited value.
Instead, we replace the right-hand side of \eqref{eq:chance-exact} by 
\begin{equation} \label{eq:chance-approx}
    \sum_{\md \in \W} \ddP_{\md}(\full{z}_k) \indi[{\Re_+}](\zeta(\md_k)) \leq \sum_{\md \in  \W} \ddP_{\md}(\full{z}_k) \sigma(\zeta(\md_k)),
\end{equation}
where 
\begin{equation}\label{eq:sigmoid-function}
    \sigma(x) = \tfrac{a}{1 + \exp\left( - \alpha (x - \overline{x} \right))}
\end{equation}
is the sigmoid function with parameters $a>0$ $\alpha>0$, $\overline{x}$.
The larger the  $\alpha$, the more accurate the approximation \eqref{eq:chance-approx} 
can be made, at the cost of larger gradients and curvature, which tend to 
impede convergence of numerical solvers. The parameters $a$ or $\overline{x}$ can be selected such that the inequality 
$\sigma \geq \indi[\Re_+]$ (and therefore \eqref{eq:chance-approx}) is satisfied uniformly by solving 
$\sigma(0) = 1$. The remaining parameter can be freely selected, e.g., to obtain the least conservative
estimate. In \cref{sec: experiment}, we simply select $a$ and $\alpha$ as tuning parameters and 
select $\overline{x}$ to satisfy \eqref{eq:chance-approx}.

We illustrate the potential benefit of using the sigmoid approximation over the 
popular $\AVAR$ approximation using a simple example. 

\begin{example}[Comparison with $\AVAR$]\label{ex:chance-constraints}
    \begin{figure}[htb!]
        \centering 
        \input{tikz/plot_avar_chance_constraint_0}
        \caption{Comparison between approximations $\nu$
        of the indicator function $\1_{\Re_+}$ in the case of \cref{ex:chance-constraints}.
        The blue histogram (left axis) represents the probability mass function of a random quantity $\zeta$; 
        The green curves (right axis) represent $\1_{\Re_+}$ and its upper bounds corresponding to $\AVAR$ and the sigmoid approximation \eqref{eq:chance-approx}.} 
        \label{fig:example-chance-constraints}
    \end{figure}
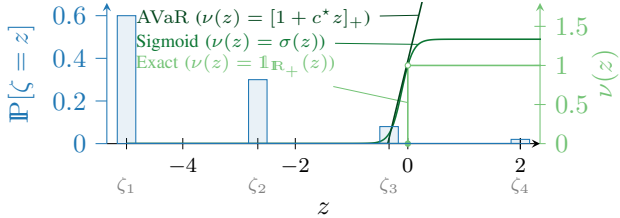
    Consider a random variable $\zeta: \Omega \to \Re$, represented by the discrete outcomes $(\zeta_i)_{i=1}^4$ (uniformly spaced between -5 and 2),
    with probabilities $(\prob\{\zeta = \zeta_i\})_{i=1}^{4} = \smallmat{0.6 & 0.3 & 0.08 & 0.02}$, as illustrated in \cref{fig:example-chance-constraints}.
    Suppose that our goal is to ensure $\prob\{\zeta > 0\} \leq \gamma = 0.05$. It is clear from the set-up that $\prob\{\zeta > 0\} = 0.02$ in this case, 
    so $\zeta$ satisfies our constraint. 

    The well-known \textit{average value-at-risk}
    is defined as 
    \(
        \AVAR_{\gamma}(\zeta) \dfn \min_{t \in \Re} t + \tfrac{1}{\gamma}\E\{[\zeta - t]_+\}, 
    \)
    with minimizer $t^\star$. Then, through the change of variables $c^\star = -\nicefrac{1}{t^\star}$, 
    it can be shown that the condition $\AVAR_{\gamma}(\zeta) \leq 0$ is equivalent to
    \(
        \E\{[1 + c^\star \zeta]_+\} \leq \gamma 
    \)
    \cite[\S 6.2.4]{shapiro_LecturesStochasticProgramming_2021}. Thus, 
    $\AVAR$ can be interpreted as an approximation of the form \eqref{eq:chance-approx}, 
    but with a piecewise-affine upper bound for the indicator $\indi[\Re_+]$.
    This is shown in \cref{fig:example-chance-constraints} for this particular example of $\zeta$.
    The compared approximations of the violation probability are thus given as
    $\sum_{i=1}^4 p_i \nu(\zeta_i)$, for $\nu \in \{[1 + c^\star (\,\cdot\,) ]_+, \sigma\}$,
    resulting in $0.14$ and $0.037$, for $\AVAR$ and $\sigma$, respectively. 
    Notice that due to the large slope $c^\star$ obtained for the $\AVAR$ approximation, the value of the largest (highly unlikely) 
    realization of $\zeta$ is heavily overweighed, resulting in an overapproximation of the constraint violation rate by a factor 7, 
    compared to a factor 1.85 for the sigmoid approximation with $\alpha=10$, $a=1.33$, $\overline{x}=-0.11$.
\end{example}

\subsection{Optimal control over scenario trees} \label{sec:scenario-tree-ocp}

\subsubsection{Scenario tree notation} \label{sec:scenario-tree-notation}
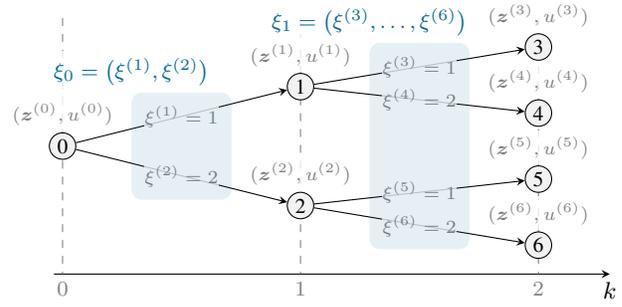
\begin{figure}[ht!]
    \centering
    \input{tikz/scenario_tree_notation}
    \caption{A fully branching scenario tree of horizon $N=2$.}
    \label{fig:scen-tree-notation}
\end{figure}

Consider an arbitrary $N$-step policy $(\pol_k)_{k=0}^{N-1} \in \Pol$. 
Since $\W = \{ 1, \dots, \nModes\}$ is a finite set, 
all realizations of the stochastic process $(\full{z}_k, u_k)_{k=1}^{N}$
satisfying \eqref{eq:full-dynamics} with $u_k = \pol_{k}\left( \full{z}_k \right)$, and 
$\full{z}_0$ known 
can be represented on a scenario tree \cite{pflug_MultistageStochasticOptimization_2014}
as illustrated in \cref{fig:scen-tree-notation}.

A scenario tree is a tree with a unique root node (with index 0) representing the current time step. 
For all possible outcomes of $\md_0$, a new node is connected to node 0. The result is a set of nodes 
representing the outcomes of the process at time $1$. Repeating this process $N$ times results in a 
tree, whose nodes are partitioned into $N+1$ time steps, or \textit{stages}. 
We denote the set of all nodes at stage $k \in \natseq{0}{N}$ by $\nodes(k)$. 
By extension, given $a < b \in \natseq{0}{N}$, 
we denote $\nodes\left(a,b\right) = \bigcup_{k=a}^{b} \nodes(k)$. 
For a given node $\iota \in \nodes(k)$, $k \in \natseq{1}{N}$, 
we denote its ancestor as $\anc{\iota} \in \nodes(k-1)$. Conversely, 
we define the children of $\iota \in \nodes\left( k \right)$, $k \in \natseq{0}{N-1}$ 
as $\ch(\iota) \dfn \{ j \in \nodes(k+1) \mid \anc(j) = \iota \}$.
A node $\iota \in \nodes(N)$ for which $\ch(\iota) = \emptyset$ is called a \textit{leaf node}.
Each leaf node $\iota$ corresponds to a scenario 
$\scen(\iota) \dfn \{\iota_k \in \nodes(k)\}_{k=0}^{N}$ 
with $\iota_{k-1}=\anc(\iota_k), \forall k \in \natseq{1}{N}$, 
which represents a single realization of the underlying stochastic process 
over $N$ future time steps. 
We denote a variable $z$ associated to a node $\iota \in \nodes\left( 0, N \right)$ 
as $\nodevar{z}{\iota}$, as illustrated in \cref{fig:scen-tree-notation}. 

\subsubsection{Optimal control problem} \label{sec:ocp-scenariotree}
Using this structure, we can represent the $N$-step policy $(\pol_{k})_{k=0}^{N-{1}} \in \Pol$ (cf. \eqref{eq: formulation_origin}) by 
a collection $\full{u} \dfn \left( \nodevar{u}{\iota}\right)_{\iota \in \nodes\left( 0, N-1 \right) } \in \Re^{\na M}$
of finite-dimensional control vectors.
Here, $M \dfn |\nodes\left( 0, N-1 \right)|$ is the number of non-leaf nodes in the scenario tree.
Accordingly, the stochastic dynamics \eqref{eq:full-dynamics} can 
be represented in scenario tree notation as 
\begin{equation} \label{eq:dynamics-tree}
    \nodevar{\full{z}}{\iota} = \fulldyn
    (
        \nodevar{\full{z}}{\anc{\iota}},
        \nodevar{u}{\anc{\iota}}, 
        \nodevar{\md}{\iota} 
    ). 
\end{equation}
Now we may decompose the control objective \eqref{eq: obj_origin} as follows. 
For a given parameter (estimate) $\param$, let 
\(
    \hat{\ddP}\left( \argdot, \param \right) 
    = \big(\hat{\ddP}_{\md}(\argdot, \param ) \big)_{\md \in \W}
\)
denote our estimate of the distribution $\ddP$, as obtained from \eqref{eq:model-probability}--\eqref{eq: param_update}. 
Then, for each leaf node $s \in \nodes(N)$, the estimated probability of its associated scenario 
is given as  
\begin{equation} \label{eq:prob-scenario-def}
    \scenprob_{s}\big( \full{u}, \theta \big) \dfn 
        \hspace{5pt}\prod_{\mathclap{\iota \in \scen\left( s \right) \setminus \{0\}}} 
            \hat{\ddP}_{\nodevar{\md}{\iota}}
                (
                \nodevar{\full{z}}{\anc{\iota}}, \param
                ), 
\end{equation}
where the predicted state $\nodevar{\full{z}}{\iota}$ in node $\iota \in \nodes\left( 0, N-1 \right)$,
can be expressed as a function of $\full{u}$ and the known 
initial state 
\(\nodevar{\full{z}}{0}\)
using \eqref{eq:dynamics-tree}. (We omitted the dependence on $\nodevar{\full{z}}{0}$ to ease notation.)
Using \eqref{eq:prob-scenario-def}, the expected cost \eqref{eq: obj_origin}
can be expressed on the scenario tree as the function $J_N: \Re^{\na M} \to \Re_+$
\[ 
    J_N\left( \full{u}, \theta \right) \dfn \> \sum_{\mathclap{s \in \nodes(N)}} \scenprob_s(\full{u}, \theta) 
    \Big( 
    \ell_N(\nodevar{\full{z}}{s}) 
        {}+{}
    \sum_{\mathclap{\substack{\iota \in \scen(s)\\ \setminus{\nodes{(N)}}}}} 
    \ell (
         \nodevar{\full{z}}{\iota}, \nodevar{u}{\iota} 
    )
    \Big).
\]

Combining these ingredients, we obtain our reformulated optimal control problem: Given $\full{\bar{z}}$, and $\theta$,
\begin{subequations}\label{eq: formulation_onestep}
    \begin{align}
        &\minimize_{\full{u} \in \Re^{\na M}} && 
            J_N (\full{u}, \theta) \label{eq: obj_onestep}\\
        &\stt&& 
        \nodevar{z}{0} = \bar{\full{z}},\\
        &&&\nodevar{\full{z}}{\iota_+} = 
            \fulldyn
            \big(
                \nodevar{\full{z}}{\iota}, \nodevar{u}{\iota}, \nodevar{\md}{\iota_+}
            \big)
            \; \forall \iota_+ \in \ch{\iota}, 
            \\
        &&&\smashoperator[l]{\sum_{\iota_+ \in \ch{\iota}}} \hat{\ddP}_{\iota_+}
            ( 
                \nodevar{z}{\iota}, \theta
            )
            \sigma_\alpha \Big(
                \overline{g}\big(
                    \nodevar{\full{z}}{\iota_+}
                \big)
            \Big) 
            \leq \gamma, \label{eq:ocp-approx-chance-constraint}
            \\
        &&&\nodevar{u}{\iota}_k \in \U, \; \nodevar{\full{z}}{\iota}_k \in \Z,\\
        &&& \nodevar{\full{z}}{\iota_{N}} \in \Z,\; \forall \iota_N \in \nodes{N},
            \\
        &&& | \nodevar{u}{\iota_+} - \nodevar{u}{\iota} | \leq \dumax\;
             \forall \iota_+ \in \ch{\iota},
            \\ \label{eq: slew_rate_constraint}
        &&&\forall \iota \in \nodes(0, N-1), 
    \end{align}
\end{subequations}
where $\overline{g}(z) \dfn \max_{k,l \in \natseq{1}{\ncircle}} g_{kl}(z)$. 
In practice, we avoid this pointwise-maximum and further approximate \eqref{eq:ocp-approx-chance-constraint}
using the union bound\footnote{
    $\prob\{ \cup_{i=1}^n E_i \} \leq \sum_{i=1}^n \prob\{E_i\}$ for any set of events 
    $\{ E_i \in \mathcal{F} \}_{i=1}^n$.}.
    In particular:
\begin{enumerate}[label=(\roman*)]
    \item If node $\iota$ has more than one child node, we approximate $\sigma \Big(
        \overline{g}\big(
            \nodevar{\full{z}}{\iota_+}
        \big)
    \Big)$
    with $\sum_{k,l \in \natseq{1}{\ncircle}}\sigma\Big(g_{kl}\big(\nodevar{\full{z}}{\iota_+}\big)\Big)$;
    \item If node $\iota$ has only one child, i.e., $\ch(\iota) = \{\iota_+\}$, 
    constraint \eqref{eq:ocp-chance-constraint} reduces to
    $g_{kl}(z^{\iota_+}) \leq 0\;\; \forall k,l \in \natseq{1}{\ncircle}$.
\end{enumerate}

This problem is a standard nonlinear problem, which can be solved with off-the-shelf numerical solvers.

\subsubsection{Scenario tree approximations}\label{sec: scenario_tree_approximation}
It is clear that number of scenarios in the tree grows exponentially with 
the prediction horizon.
Several techniques exist to combat this explosive growth while 
maintaining a sufficiently long prediction horizon for planning 
complex tasks. For instance, in \cite{bernardini_StabilizingModelPredictive_2012}, 
a scenario tree is dynamically built by repeatedly adding the most 
likely child node to any of the current leaf nodes, until 
a desired number of nodes is reached.
However, this approach is not directly applicable in our setting, 
since the transition probabilities are not yet available when 
constructing the tree,
as they depend on the future states over which we need to optimize.

\begin{figure}[ht!]
    \centering 
    \begin{minipage}{0.32\columnwidth}
    \begin{tikzpicture}
        \scenariotree[horizon=8,branching_level=2,hspace=5pt,vspace=4pt]{}
    \end{tikzpicture}
    {\scriptsize $(N, \bhor, \td) = (8, 2,1)$}
    \end{minipage}
    \hfill
    \begin{minipage}{0.32\columnwidth}
        \begin{tikzpicture}
            \scenariotree[horizon=8,branching_level=4, skips=2,hspace=5pt,vspace=4pt]{}
        \end{tikzpicture}
    {\scriptsize $(N, \bhor, \td) = (8, 4, 2)$}
    \end{minipage}
    \hfill
    \begin{minipage}{0.32\columnwidth}
        \begin{tikzpicture}
            \scenariotree[horizon=8,branching_level=8, skips=4,hspace=5pt,vspace=4pt]{}
        \end{tikzpicture}
        {\scriptsize $(N, \bhor, \td) = (8, 8, 4)$}
    \end{minipage}
    \caption{Different combinations of the branching horizon $\bhor$ and timescale factor $\td$ 
    leading to an approximate scenario tree of horizon $N$, each with the same number of scenarios.}
    \label{fig:scenariotrees-approx}
\end{figure}

More direct approaches (illustrated in \cref{fig:scenariotrees-approx}) involve a combination of the following:
\begin{inlinelist*}
    \item Fix a \textit{branching horizon} $\bhor < N$, after which the modes are kept constant, i.e., 
    set $\md_{k} = \md_{k-1}$ for all $k > \bhor$ \cite{schuurmans_SafeLearningBasedMPC_2022b}
    \item assume that the discrete process $(\md_t)_{t\in\N}$ governing the 
    branching, evolves at a slower timescale than the continuous state updates. 
    That is, assume that $\md_{k} \sim \ddP\left( \full{z}_k \right)$ if $k\mod\td = 0$, otherwise, 
    $\md_{k}=\md_{k-1}$ \cite{chen_InteractiveMultiModalMotion_2022}.
\end{inlinelist*}
In general, the most suitable choice depends on the particular application.  
We have found that for our use case, the resulting behavior is very similar in both cases.

%% file: tikz/plot_avar_chance_constraint_0.tex
\usepgfplotslibrary{groupplots}

\begin{tikzpicture}

\definecolor{darkgray176}{RGB}{176,176,176}
\definecolor{darkgreen06827}{RGB}{0,68,27}
\definecolor{darkseagreen115195117}{RGB}{115,195,117}
\definecolor{forestgreen1111851}{RGB}{11,118,51}
\definecolor{seagreen5516085}{RGB}{55,160,85}
\definecolor{steelblue31119180}{RGB}{31,119,180}
\newlength{\avarfigheight}
\newlength{\avarfigwidth}
\setlength{\avarfigheight}{0.4\columnwidth}
\setlength{\avarfigwidth}{0.85\columnwidth}

\tikzset{point/.style={circle, fill, minimum size=3pt, inner sep=0pt}}

\begin{axis}[
ylabel shift=-0.2cm,
xlabel shift=-0.2cm,
tick align=outside,
axis y line=right,
axis x line=bottom,
height=\avarfigheight,
width=\avarfigwidth, 
x grid style={darkgray176},
xlabel={$z$},
xmin=-5.35, xmax=2.35,
xtick style={color=black},
y grid style={darkseagreen115195117},
y axis line style={color=darkseagreen115195117},
x tick label style={font=\footnotesize},
y tick label style={color=darkseagreen115195117},
x tick style={color=black},
ylabel=\textcolor{darkseagreen115195117}{\(\nu(z)\)},
extra x ticks={-5.0, -2.66666667, -0.33333333,2.0},
extra x tick style={
    tick label style={
        yshift=-0.25cm,
        opacity=0.5,
    }
},
extra x tick label={
    \scriptsize
    \pgfmathparse{int(\ticknum+1)}
    $\zeta_\pgfmathresult$
},
ymin=-0.0, ymax=1.8,
ytick style={color=darkseagreen115195117}
]
\addplot [semithick, darkseagreen115195117]
table {%
-5 0
0 0 
0 1 
2.5 1
}
coordinate[pos=0.65] (coord_exact)
node[anchor=west] (label_exact) at (axis cs: -5, 1.) {\scriptsize Exact ($\nu(z) = \indi[\Re_+](z)$)
};
\draw[very thin, darkseagreen115195117] (coord_exact) -- (label_exact);

\addplot [only marks, mark size=1, darkseagreen115195117]
table {%
0 0
};
\addplot [only marks, fill=white, mark size=1, draw=darkseagreen115195117]
table {%
0 1
};

\addplot [semithick, forestgreen1111851]
table {%
-5 7.68142200916094e-22
-4.96482412060301 1.09196452271497e-21
-4.92964824120603 1.55229919336036e-21
-4.89447236180905 2.20669512203214e-21
-4.85929648241206 3.13696185788713e-21
-4.82412060301507 4.45939704111755e-21
-4.78894472361809 6.33932539547102e-21
-4.75376884422111 9.01176685976206e-21
-4.71859296482412 1.28108176924829e-20
-4.68341708542714 1.8211417639178e-20
-4.64824120603015 2.58887247004673e-20
-4.61306532663317 3.68025202593091e-20
-4.57788944723618 5.23171965057204e-20
-4.5427135678392 7.43723264312547e-20
-4.50753768844221 1.05725140264201e-19
-4.47236180904523 1.50295221626784e-19
-4.43718592964824 2.13654515731987e-19
-4.40201005025126 3.03723908176032e-19
-4.36683417085427 4.31763457382019e-19
-4.33165829145729 6.13780074970013e-19
-4.2964824120603 8.72528635736026e-19
-4.26130653266332 1.24035668674412e-18
-4.22613065326633 1.76324838789166e-18
-4.19095477386935 2.50657323867348e-18
-4.15577889447236 3.56325827034878e-18
-4.12060301507538 5.06540535313797e-18
-4.08542713567839 7.20080596040736e-18
-4.05025125628141 1.02364179891975e-17
-4.01507537688442 1.45517396004987e-17
-3.97989949748744 2.06862523222657e-17
-3.94472361809045 2.94068645322503e-17
-3.90954773869347 4.18037867926095e-17
-3.87437185929648 5.94268249267248e-17
-3.8391959798995 8.44791295676573e-17
-3.80402010050251 1.20092623849733e-16
-3.76884422110553 1.70719541938023e-16
-3.73366834170854 2.42689026729874e-16
-3.69849246231156 3.44998369996069e-16
-3.66331658291457 4.90437812140655e-16
-3.62814070351759 6.97189518837592e-16
-3.5929648241206 9.91100631200086e-16
-3.55778894472362 1.4089145556906e-15
-3.52261306532663 2.00286445467523e-15
-3.48743718592965 2.84720319454377e-15
-3.45226130653266 4.04748609527587e-15
-3.41708542713568 5.75376696782496e-15
-3.38190954773869 8.17935714681601e-15
-3.34673366834171 1.16274926859717e-14
-3.31155778894472 1.65292435255692e-14
-3.27638190954774 2.34974038605242e-14
-3.24120603015075 3.34031008334099e-14
-3.20603015075377 4.74846988165126e-14
-3.17085427135678 6.75026139920409e-14
-3.1356783919598 9.59593934325205e-14
-3.10050251256281 1.36412571949034e-13
-3.06532663316583 1.9391941862196e-13
-3.03015075376884 2.75669173166294e-13
-2.99497487437186 3.9188181139473e-13
-2.95979899497487 5.57085699275376e-13
-2.92462311557789 7.91933861978904e-13
-2.8894472361809 1.1257859294621e-12
-2.85427135678392 1.60037854147028e-12
-2.81909547738693 2.27504306899793e-12
-2.78391959798995 3.23412294758673e-12
-2.74874371859296 4.59751790312742e-12
-2.71356783919598 6.53567326045545e-12
-2.67839195979899 9.29088822870055e-12
-2.64321608040201 1.3207607026566e-11
-2.60804020100502 1.87754797037718e-11
-2.57286432160804 2.66905759231832e-11
-2.53768844221106 3.79424043672405e-11
-2.50251256281407 5.39376165319697e-11
-2.46733668341709 7.66758597841254e-11
-2.4321608040201 1.08999764016416e-10
-2.39698492462312 1.54950314074899e-10
-2.36180904522613 2.20272034973814e-10
-2.32663316582915 3.13131145813953e-10
-2.29145728643216 4.45136462681429e-10
-2.25628140703518 6.32790679078624e-10
-2.22110552763819 8.99553456243444e-10
-2.18592964824121 1.27877424134882e-09
-2.15075376884422 1.81786146091905e-09
-2.11557788944724 2.58420930268868e-09
-2.08040201005025 3.67362302466699e-09
-2.04522613065327 5.22229608440794e-09
-2.01005025125628 7.42383641347099e-09
-1.9748743718593 1.05534703825432e-08
-1.93969849246231 1.50024503238073e-08
-1.90452261306533 2.13269670775281e-08
-1.86934673366834 3.03176823785905e-08
-1.83417085427136 4.30985736911733e-08
-1.79899497487437 6.12674485972596e-08
-1.76381909547739 8.70956951611323e-08
-1.7286432160804 1.23812240062569e-07
-1.69346733668342 1.76007213482282e-07
-1.65829145728643 2.50205784756483e-07
-1.62311557788945 3.55683906474691e-07
-1.58793969849246 5.05627946353434e-07
-1.55276381909548 7.18783176232978e-07
-1.51758793969849 1.02179719968646e-06
-1.48241206030151 1.45255127266245e-06
-1.44723618090452 2.06489596313686e-06
-1.41206030150754 2.93538313898127e-06
-1.37688442211055 4.17283580307334e-06
-1.34170854271357 5.9319519845424e-06
-1.30653266331658 8.43264304299015e-06
-1.2713567839196 1.198752328658e-05
-1.23618090452261 1.70409861458142e-05
-1.20100502512563 2.42247492306783e-05
-1.16582914572864 3.44368066097778e-05
-1.13065326633166 4.8953647148503e-05
-1.09547738693467 6.95897448651826e-05
-1.06030150753769 9.89242134614095e-05
-1.0251256281407 0.000140622866196515
-0.989949748743719 0.000199895753377759
-0.954773869346734 0.000284146992565857
-0.919597989949748 0.000403897354504971
-0.884422110552764 0.000574093375928276
-0.849246231155779 0.000815963497318776
-0.814070351758794 0.00115964689481301
-0.778894472361809 0.00164791079792045
-0.743718592964824 0.00234139513704813
-0.708542713567839 0.00332598816059304
-0.673366834170855 0.00472315074036241
-0.63819095477387 0.00670427142124121
-0.603015075376884 0.00951043151773172
-0.567839195979899 0.013479229858924
-0.532663316582915 0.0190804112577413
-0.49748743718593 0.0269616407289276
-0.462311557788945 0.0380042428003005
-0.42713567839196 0.0533849866888907
-0.391959798994975 0.0746323791455711
-0.35678391959799 0.103652515070473
-0.321608040201005 0.142679726461504
-0.28643216080402 0.194085816421857
-0.251256281407035 0.259975581463991
-0.21608040201005 0.341539527883172
-0.180904522613066 0.438263485366669
-0.14572864321608 0.547293625524959
-0.110552763819095 0.663388263893546
-0.075376884422111 0.77974057183066
-0.0402010050251258 0.889483980643692
-0.00502512562814061 0.987225172658286
0.0301507537688446 1.06992940203284
0.0653266331658289 1.13692994148703
0.100502512562814 1.18932077137446
0.135678391959799 1.22916494621159
0.170854271356784 1.25883146029279
0.206030150753769 1.2805731778566
0.241206030150754 1.2963228873393
0.276381909547738 1.30763615764567
0.311557788944723 1.31571351706644
0.346733668341709 1.32145558272619
0.381909547738694 1.32552495877169
0.417085427135678 1.32840260943568
0.452261306532663 1.33043438962025
0.487437185929648 1.33186737274886
0.522613065326633 1.33287725530102
0.557788944723618 1.33358857502608
0.592964824120603 1.33408940778486
0.628140703517587 1.33444194392285
0.663316582914573 1.33469004704254
0.698492462311558 1.334864630409
0.733668341708543 1.33498746840639
0.768844221105527 1.3350738923146
0.804020100502512 1.33513469389823
0.839195979899498 1.33517746807692
0.874371859296482 1.33520755920555
0.909547738693467 1.33522872761596
0.944723618090452 1.33524361893172
0.979899497487438 1.33525409442349
1.01507537688442 1.33526146350433
1.05025125628141 1.33526664733134
1.08542713567839 1.33527029391742
1.12060301507538 1.33527285911952
1.15577889447236 1.33527466361627
1.19095477386935 1.33527593299203
1.22613065326633 1.33527682593557
1.26130653266332 1.33527745407723
1.2964824120603 1.33527789594369
1.33165829145729 1.33527820677472
1.36683417085427 1.33527842542884
1.40201005025126 1.33527857924107
1.43718592964824 1.33527868744027
1.47236180904523 1.33527876355298
1.50753768844221 1.33527881709446
1.5427135678392 1.3352788547582
1.57788944723618 1.33527888125275
1.61306532663317 1.33527889989034
1.64824120603015 1.33527891300095
1.68341708542714 1.3352789222236
1.71859296482412 1.33527892871127
1.75376884422111 1.33527893327503
1.78894472361809 1.3352789364854
1.82412060301507 1.33527893874373
1.85929648241206 1.33527894033236
1.89447236180905 1.33527894144987
1.92964824120603 1.33527894223599
1.96482412060301 1.33527894278899
2 1.33527894317799
2.5 1.33528
}
coordinate[pos=0.72] (safe_coord)
node[anchor=west] (label_safe) at (axis cs: -5, 1.3)  {\scriptsize Sigmoid ($\nu(z) = \sigma(z)$)}
;
\draw[very thin, forestgreen1111851] (safe_coord) -- (label_safe); 
\addplot [semithick, darkgreen06827]
table {%
-5 0
-0.35678391959799 0
2 6.99999846505605
}
coordinate[pos=0.53] (avar_coord)
node[anchor=west] (label_avar) at (axis cs: -5, 1.6) {\scriptsize$\AVAR$ ($\nu(z) = [1 + c^\star z]_+$)};
\draw[very thin, darkgreen06827] (avar_coord) -- (label_avar);
\end{axis}

\begin{axis}[
axis y line=left,
axis x line=none,
tick align=outside,
height=\avarfigheight,
width=\avarfigwidth,
x grid style={darkgray176},
xmin=-5.35, xmax=2.35,
x axis line style=none,
y grid style={darkgray176},
y tick style={color=steelblue31119180},
y axis line style={color=steelblue31119180},
y tick label style={color=steelblue31119180},
ylabel=\textcolor{steelblue31119180}{\(\prob[\zeta = z]\)},
ymin=0, ymax=0.66,
ytick pos=left,
yticklabel style={anchor=east}
]







\addplot[ybar, draw=steelblue31119180, fill=steelblue31119180, fill opacity=0.1, draw opacity=1, bar width=7pt] 
table {
    -5 0.6
    -2.66666666666667 0.3
    -0.333333333333333 0.08
    2 0.02
};

\addplot [semithick, steelblue31119180]
table {%
-5 0
2 0
};
\end{axis}

\end{tikzpicture}

%% file: tikz/scenario_tree_notation.tex
\usetikzlibrary{decorations.pathreplacing}
\colorlet{modecolor}{MidnightBlue}
\begin{tikzpicture}
    {\footnotesize
    \tikzset{
        label/.style={draw=none, font=\scriptsize, inner sep=0, outer sep=0, node distance=20pt}, 
        toplabel/.style={draw=none, 
                     font=\scriptsize,
                     inner sep=0, outer sep=5pt,
                     node distance=13pt, anchor=south, color=gray,
                     fill=white},
        braces/.style={decorate,decoration={brace,amplitude=4pt},xshift=4pt, gray}, 
        braceslabel/.style={align=center, anchor=west, gray, font=\scriptsize, midway, inner sep=10pt},
        conn/.style={fill=white, fill opacity=0.8, text opacity=1,
                    inner sep=2pt, draw=none, text=gray}
    };
    
    
    \scenariotreecoordinates[%
        hspace=10em, vspace=1pt, nodesize=10pt, vspace=15pt,
        branching_factor=2, branching_level=2,horizon=2,skips=1
    ]{}


    \drawtimeline{$k$}{1.15}
    \foreach \t in {0,...,\hor}{
        \annotatestage{\t}{$\t$}{}
    }

    \fillsubtree{0}{\hor}{\name}
    
    \newcounter{nodecnt}
    \setcounter{nodecnt}{0}
    \newcommand\annotatenode[1]{
        \node (labelnode/#1) at (#1) {\thenodecnt};
        \node[above of=labelnode/#1, label, gray, conn, node distance=12pt] {$(\nodevar{\full{z}}{\thenodecnt}, \nodevar{u}{\thenodecnt})$};
        \stepcounter{nodecnt}
    }
    \doforeverynode{\annotatenode}{0}{\hor}{}

    \newcounter{modenode}\setcounter{modenode}{1}
    \providecommand\labelconnection[4]{
        \foreach \lcltimectr in {#1,...,0}
        {
            \computenchildren{\lcltimectr} 
            \ifnum \nc>1
                \breakforeach 
            \fi
        }
        \path (node/tree/#1#2) --
              (node/tree/#3#4)
              node[conn,label,midway] (connect\themodenode) {
                \pgfmathparse{int(\nc - mod(#4, \nc))}  
                $\nodevar{\md}{\themodenode}=\pgfmathresult$
            };
        \stepcounter{modenode}
    }
    \pgfmathparse{int(\hor-1)}
    \doforeveryconnection{\labelconnection}{0}{\pgfmathresult}

    \fill[opacity=0.1, modecolor, rounded corners] ($(connect1.north west) - (5pt,-5pt)$) coordinate(topleft) rectangle ($(connect2.south east) + (5pt,-5pt)$); 
    \node[text=modecolor, anchor=south, fill opacity=0.8,fill=white, text opacity=1] at (topleft) {
        $\md_0=\left(\nodevar{\md}{1},\nodevar{\md}{2}\right) $
    };
    
    \fill[opacity=0.1, modecolor, rounded corners] ($(connect3.north west) - (5pt,-5pt)$) coordinate(topleft) rectangle ($(connect6.south east) + (5pt,-5pt)$); 
    \node[text=modecolor, anchor=south, fill opacity=0.8,fill=white, text opacity=1] at (topleft) {
        $\md_1 = \left(\nodevar{\md}{3}, \dots, \nodevar{\md}{6}\right) $};

    }
\end{tikzpicture}

%% file: content/experimental_results.tex
We illustrate the proposed control method on
a simple lane-change scenario, using 
the \texttt{highway-env} \cite{highway-env} simulator.

\subsection{Interactive Target Vehicle Behavior}\label{sec: tv_model}
\begin{table}[ht!]
    \centering
    \caption{Target vehicle parameters}
    \label{tab: tv_params}
    {\footnotesize
    \begin{tabular}{lll}\toprule
    Symbol & Description & Value \\
        \midrule
    $\amaxcomf$  & maximum comfort acceleration & 3 [\mpss{}]\\
    $\kbrake, \; \ktrack$ & control gain in maneuver & 0.7 [-]\\
    $\Np$   &P-IDM prediction horizon & $U[0.1, 1]$ [\SI{}{\second}]\\
    $\cthres$     & Decision threshold & $U[0, 4]$ [\SI{}{\meter}]\\
    \bottomrule
    \end{tabular}
    }
\end{table}
\subsubsection{Maneuvers}
For simplicity, we will assume the target vehicle only moves in the longitudinal direction.
Lateral motion is handled internally by the simulator.
We consider two base policies, represented by the set $\W = \{\mdbrake, \mdtrack\}$.
That is, we assume that the (stochastic) target vehicle policy is given as
$\idx{\kappa}{\tv}(\full{z}, \md) = \smallmat{
    \max\{
        \amin,
        \min\{
            \amaxcomf, \idx{a}{\tv}(\full{z}, \md)
        \}
    \} 
        & 0}^\top$
with\footnote{See \cref{tab: tv_params} for the values used to describe the target vehicle behavior.}
\[
\idx{a}{\tv}(\full{z}, \md) =
\begin{cases}
    -k_b \idx{v}{\tv} & \text{if} \quad \md = \mdbrake \\
    k_t(v_{\max} - \idx{v}{\tv}) & \text{if} \quad \md = \mdtrack.
\end{cases}
\]
\subsubsection{Interaction model}
Following the \textit{Predictive Intelligent Driver Model} (P-IDM) \cite{brito2022learning}, the target vehicle will predict the ego vehicle's trajectory with constant velocity to make decisions.
If \begin{inlinelist*}
    \item $\pos_x^\ego > \pos_x^\tv$
    \item $\exists k\in \natseq{0}{\Np}: \; \big| \pos_{y, k}^\ego - \pos_{y, k}^\tv \big| \leq  \cthres$
\end{inlinelist*}
the target vehicle will choose \mdbrake. Otherwise it chooses \mdtrack.
The prediction horizon $\Np$ and the threshold $\cthres$ may differ per individual driver.
Note that this behavior can be captured quite well using the proposed probabilistic model \eqref{eq:model-probability}.

\subsection{Simulation Results and Analysis}\label{sec: experiment_results}
We consider a straight two-lane highway. 
Each vehicle is modeled as a box with a length of \SI{5}{\meter} and a width of \SI{2}{\meter}.

\subsubsection{Simulation overview}
\begin{figure}[tb]
    \centering
    \includegraphics[width=\linewidth]{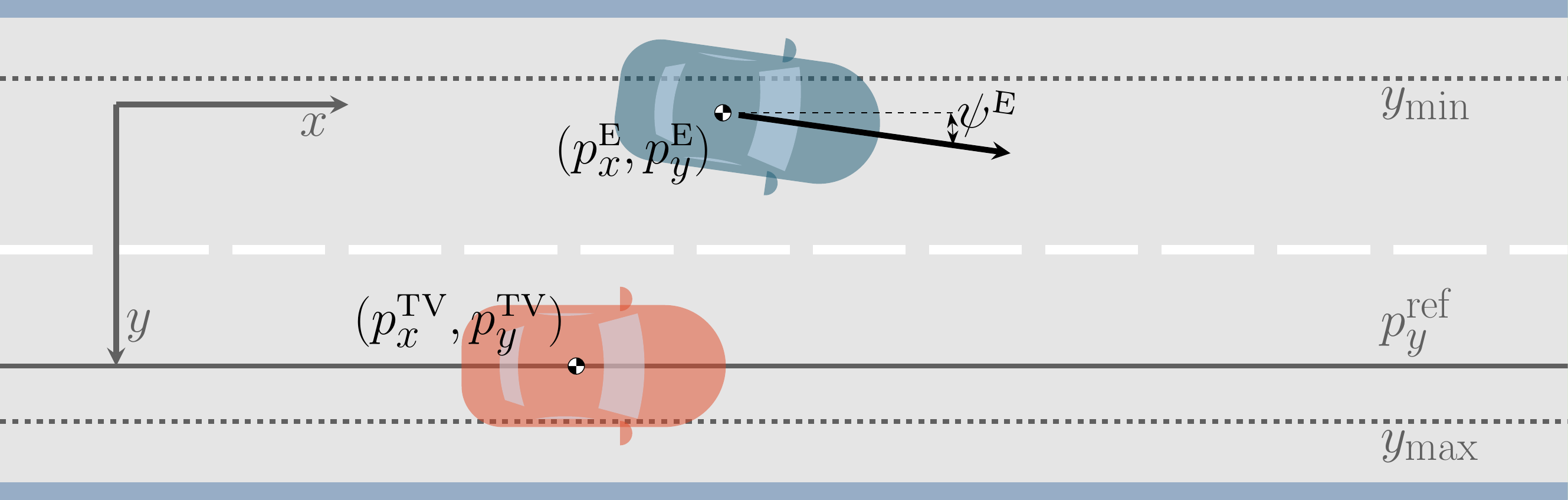}
    \caption{Illustration of the coordinate system.}
    \label{fig: illustration_setup}
\end{figure}
The goal of ego vehicle is to drive in the center of the right lane at the maximum velocity, i.e., 
the reference state in \eqref{eq:stage-cost} is
\(
z^{\refer}=(0, \pos_y^{\refer}, v^{\refer},0) 
\)
(cf. \cref{tab: mpc_params}). 
The coordinate system is illustrated in \cref{fig: illustration_setup}.
The two vehicles are initialized on two adjacent lanes with random relative distance $\pos_{x,\init}^\ego - \pos_{x,\init}^\tv \sim U[0, 5]\SI{}{\meter}$, $\pos_{y,\init} = \pos_{y, \init}^\ego - \pos_{y,\init}^\tv \sim U[-5, -3]\SI{}{\meter}$. 
Over all experiments, the ego vehicle has fixed longitudinal initial position $\pos_{x, \init}^\ego = \SI{6}{\meter}$ and the target vehicle has fixed lateral initial position $\pos_{y, \init}^\tv = \SI{4}{\meter}$. 
The initial velocities are drawn randomly as $v_{\init}^\ego, v_{\init}^\tv \sim U[23, 25]\mps{}$. 
The hyperparameters and constraints used in the simulation are summarized in \cref{tab: mpc_params}.
As described in \S\ref{sec: scenario_tree_approximation}, we control the scenario tree complexity by selecting branching horizon $\bhor = 11$ and timescale factor $D=5$.
\begin{table}[tb]
    \centering
    \caption{\ac{MPC} parameters}
    \label{tab: mpc_params}
    {\footnotesize
    \begin{tabular}{lll}\toprule
    Symbol & Description & Value \\
        \midrule
    $\Ts$                           & sampling time & \SI{0.1}{\second} \\
    $T_{\mathrm{sim}}$                           & simulation horizon & \SI{6}{\second} \\
    $Q $                              &  weights for states & $\diag{0, 1, 0.01, \nicefrac{16}{\pi^2}}$\\
    $R $                              &  weights for input  & $\diag{0.01, \nicefrac{16}{\pi^2}}$\\
    $N $                              & prediction horizon & 20 \\
    $r_c$                               &  circle radius in \eqref{eq: collision_avoidance_robust} & \SI{1.3}{\meter}\\
    $(\alpha, a)$                      & parameters in \eqref{eq:sigmoid-function}  & (10, 1.2)\\
    $\gamma$                        & threshold in \eqref{eq:ocp-approx-chance-constraint} & 0.05\\
    $L$    & window length in \eqref{eq: param_update} & 15 \\
    $\lambda$ & weight in \eqref{eq: param_update} & 1 \\
    $[y_{\min}, y_{\max}]$ &bounds on lateral position & $[-1, 5]\SI{}{\meter}$\\
    $[\vmin, \vmax]$ & bounds on velocity & $[0, 28]\mps{}$ \\
    $[\psi_{\min}, \psi_{\max}]$ & bounds on heading angle &$[-\nicefrac{\pi}{4}, \nicefrac{\pi}{4}]\SI{}{\radian}$ \\
    $[\amin, \amax]$ & bounds on acceleration & $[-5, 5]\mpss{}$ \\
    $[\delta_{\min}, \delta_{\max}]$ & bounds on steering angle &$[-\nicefrac{\pi}{4}, \nicefrac{\pi}{4}]\SI{}{\radian}$ \\ 
    $\epsilon_{\delta_u}$  & bounds on slew rate \eqref{eq: slew_rate_constraint} & $[5, \nicefrac{\pi}{4}]$ \\
    \bottomrule
    \end{tabular}
    }
\end{table}

We compare the following variants of the method:
\begin{enumerate}[label=(\roman*)]
    \item \MLE: Estimate $\param_t$ according to \eqref{eq: param_update}, with $\param_0 = 0$;
    \item \MAP: Estimate $\param_t$ according to \eqref{eq: param_update}, with $\param_0 = \hat{\param}$ obtained according to \cref{sec: offline_estimation};
    \item \PRIOR: Use the prior $\hat{\param}$ without online learning;
    \item \EMP: Replace the state-dependent distribution by the empirical distribution 
    \(
        \hat{P}_t = \big( t^{-1} \sum_{j=0}^{t-1} \indi[\{i\}](\md_j) \big)_{i \in \W}
    \) 
    with $t$ the current time step of the simulation;
    \item \UNI: Use uniform distribution $\hat{P}_t = (0.5, 0.5)$, $\forall t \in \N$;
    \item \BRA: Assume $\md_{t} = \mdbrake$, i.e., $\hat{P}_t = (1, 0)$, $\forall t \in \N$;
    \item \TRA: Assume $\md_{t} = \mdtrack$, i.e., $\hat{P}_t = (0, 1)$, $\forall t \in \N$.
\end{enumerate}
Each method is simulated 50 times with random initial states and target vehicle parameters.
The problem \eqref{eq: formulation_onestep} is implemented with \casadi{} \cite{Andersson2019} and solved using \IPOPT{} \cite{wachter2006implementation}.
After the first step, the solver is warm-started with a shifted version of the previous solution. 
If the problem is infeasible with warm-starting, it is resolved with zero as initial guess.
The convex parameter estimation problem \eqref{eq: param_update} is solved using \MOSEK{}\cite{mosek} through \CVXPY{} \cite{diamond2016cvxpy}.

\subsubsection{Offline parameter estimation}\label{sec: offline_estimation}
An informative initial guess $\hat{\param}$ for the parameter is estimated from a synthetic offline dataset 
containing 800 training points and 200 validation points generated from 10 different drivers.
We solve \eqref{eq:offline-param} to obtain \eqref{eq:model-probability} with feature function $\feat(\full{z}) = z^\ego - z^\tv$.
The fitted function achieves a misclassification rate of 0.175 on the validation set.  
We point out that given the same states, different drivers may choose different maneuvers according to their own driving style.
Therefore a model with zero prediction error is not expected.

\subsubsection{Results and analysis}
\begin{figure}[tb]
    \centering
    \input{tikz/simulation_illustration}

    \vspace{-5pt}
    \caption{Snapshot of a closed-loop simulation.
    The ego vehicle and target vehicle are shown in blue and yellow, respectively.
    The last 5 time steps are shown as transparent rectangles and the predicted trajectories are shown as lines}
    \label{fig: snapshot}
\end{figure}
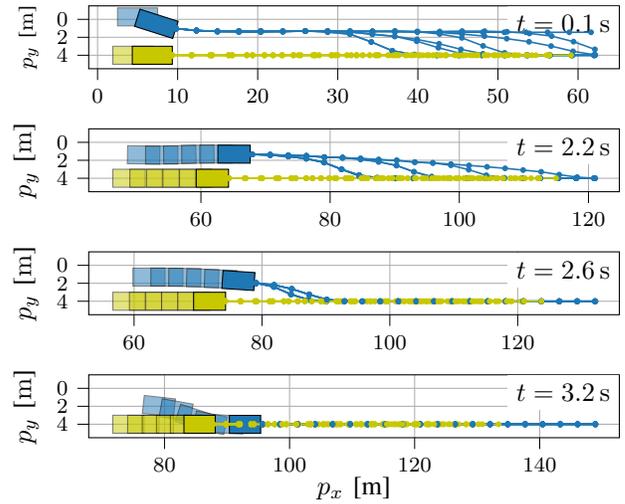

A snapshot of successful lane-changing is presented in \cref{fig: snapshot}. 
The control performance is quantified with closed-loop cost. 
The simulation result presented in \cref{tab: result} demonstrates the advantage of employing the decision-dependent distribution instead of a fixed distribution, and furthermore, \MAP{} estimator performs the best on average.
This shows the benefit of obtaining an informative prior estimate.
On the other hand, without informative prior, online learning can efficiently capture the driver behavior and reaches better level control performance than only using offline approximated distribution.
In addition, no collision occurs in the simulations, which proves the effectiveness of our collision avoidance constraint formulation.
\begin{table}[ht!]
    \centering
    \caption{Simulation results.
    If the ego vehicle is sufficiently close to the reference$^{1}$, we record the relative positive of ego vehicle to the target vehicle (front / behind).
    Otherwise we categorize the run under `time-out'.}
    \label{tab: result}
    \begin{tabular}{lrrrrrr}\toprule
        & \multicolumn{4}{c}{number of simulations} & \multicolumn{2}{c}{closed-loop cost}
\\\cmidrule(lr){2-5}\cmidrule(lr){6-7}
      method &  collision &  front &  behind &  time-out &      mean &      Q3$^{\hyperlink{Q3}{2}}$ \\
     \midrule
         \MLE &         0 &         34 &           5 &            11 & 154.31 & 163.30	 \\
         \MAP &         0 &         31 &          11 &            8 &  139.50 &  161.77	\\
         \PRIOR &         0 &         25 &          14 &            11 & 157.13 & 165.11\\
         \EMP &         0 &         34 &           5 &            11 & 154.85 & 160.26\\
         \UNI &         0 &         33 &           5 &            12 & 167.89 & 243.92 \\
         \BRA &         0 &         33 &           1 &            16 & 186.86 & 302.64\\
         \TRA &         0 &         25 &          16 &             9 & 159.54 & 211.21 \\
     \bottomrule
    \end{tabular}
    
    \vspace{2pt}
    \raggedleft{
    {\footnotesize
    \hypertarget{close}{
        $^1\big( 
        \abs{\pos_y^\ego - \pos_y^\refer} \leq 0.1$ and $\abs{\head^\ego - \head^\refer} \leq 0.01
        \big)$}
        \qquad \hypertarget{Q3}{$^2$\textit{Third quartile}}
    }
    }
\end{table}

%% file: tikz/simulation_illustration.tex
\begin{tikzpicture}

\definecolor{darkgray176}{RGB}{176,176,176}
\definecolor{goldenrod2002000}{RGB}{200,200,0}
\definecolor{steelblue31119180}{RGB}{31,119,180}
\newlength{\simulationheight}
\newlength{\simulationwidth}
\newlength{\markersize}
\setlength{\simulationheight}{0.28\columnwidth}
\setlength{\simulationwidth}{1\columnwidth}
\setlength{\markersize}{0.8pt}

\begin{groupplot}[group style={group size=1 by 4, vertical sep=8mm}]
\nextgroupplot[
height=\simulationheight,
width=\simulationwidth,
tick align=outside,
tick pos=left,
x grid style={darkgray176},
xmajorgrids,
xmin=-1.10173257443619, xmax=65.1587567857857,
xtick style={color=black},
y dir=reverse,
y grid style={darkgray176},
ylabel={\(\displaystyle p_y\) [$\SI{}{\meter}$]},
label style={font=\footnotesize},
ticklabel style={font=\footnotesize},
ymajorgrids,
ymin=-1.5, ymax=5.5,
ytick style={color=black}
]
\node[fill=white, anchor=north east] at (rel axis cs:1,1) {$t = \SI{0.1}{\second}$};
\draw[draw=black,fill=steelblue31119180,opacity=0.5] (axis cs:2.5,-1.27096974849701) rectangle (axis cs:7.5,0.729030251502991);
\draw[draw=black,fill=steelblue31119180,rotate around={-17.6665694486358:(axis cs:7.34422636032104,0.49987980723381)}] (axis cs:4.84422636032104,-0.50012019276619) rectangle (axis cs:9.84422636032104,1.49987980723381);
\draw[draw=black,fill=goldenrod2002000,opacity=0.5] (axis cs:1.91010785102844,3.00000011920929) rectangle (axis cs:6.91010785102844,5.00000011920929);
\draw[draw=black,fill=goldenrod2002000] (axis cs:4.34008938074112,2.99999997019768) rectangle (axis cs:9.34008938074112,4.99999997019768);
\addplot [semithick, steelblue31119180, mark=*, mark size=\markersize, mark options={solid}]
table {%
7.34422636032104 0.49988
9.80523636032104 1.031227
12.363506360321 1.25123
14.979466360321 1.346801
17.646976360321 1.313773
20.364686360321 1.3104
23.131906360321 1.36288
25.930856360321 1.286068
28.727666360321 1.419509
31.518366360321 1.647573
34.180096360321 2.516595
36.802236360321 3.498628
39.568906360321 3.929349
42.367536360321 4.015752
45.167476360321 4.011721
47.967376360321 4.002956
50.767266360321 4.000086
53.567146360321 3.99978
56.367016360321 3.999912
59.166876360321 3.999986
61.9667363603211 4.000006
};
\addplot [semithick, goldenrod2002000, mark=*, mark size=\markersize, mark options={solid}]
table {%
6.84008636032105 4
9.30006636032105 4
11.710056360321 4
14.070036360321 4
16.380016360321 4
18.639996360321 4
20.849976360321 4
23.009956360321 4
25.119946360321 4
27.179926360321 4
29.189906360321 4
31.149886360321 4
33.059866360321 4
34.919846360321 4
36.729826360321 4
38.489816360321 4
40.199796360321 4
41.859776360321 4
43.469756360321 4
45.029736360321 4
46.539716360321 4
};
\addplot [semithick, steelblue31119180, mark=*, mark size=\markersize, mark options={solid}]
table {%
7.34422636032104 0.49988
9.80523636032104 1.031227
12.363506360321 1.25123
14.979466360321 1.346801
17.646976360321 1.313773
20.364686360321 1.3104
23.131906360321 1.36288
25.930856360321 1.286068
28.727666360321 1.419509
31.518366360321 1.647573
34.180096360321 2.516595
36.802236360321 3.498628
39.568906360321 3.929349
42.367536360321 4.015752
45.167476360321 4.011721
47.967376360321 4.002956
50.767266360321 4.000086
53.567146360321 3.99978
56.367016360321 3.999912
59.166876360321 3.999986
61.9667363603211 4.000006
};
\addplot [semithick, goldenrod2002000, mark=*, mark size=\markersize, mark options={solid}]
table {%
6.84008636032105 4
9.30006636032105 4
11.710056360321 4
14.070036360321 4
16.380016360321 4
18.639996360321 4
20.849976360321 4
23.009956360321 4
25.119946360321 4
27.179926360321 4
29.189906360321 4
31.149886360321 4
33.139866360321 4
35.159846360321 4
37.209826360321 4
39.289816360321 4
41.399796360321 4
43.539776360321 4
45.709756360321 4
47.909736360321 4
50.139716360321 4
};
\addplot [semithick, steelblue31119180, mark=*, mark size=\markersize, mark options={solid}]
table {%
7.34422636032104 0.49988
9.80523636032104 1.031227
12.363506360321 1.25123
14.979466360321 1.346801
17.646976360321 1.313773
20.364686360321 1.3104
23.131906360321 1.36288
25.931316360321 1.305675
28.731176360321 1.333971
31.529986360321 1.415505
34.3191263603211 1.661905
37.109316360321 1.896138
39.856416360321 2.437826
42.560206360321 3.16549
45.290276360321 3.787356
48.082856360321 3.991044
50.882716360321 4.014583
53.682636360321 4.006102
56.482546360321 4.000996
59.282446360321 3.999839
62.082346360321 3.999827
};
\addplot [semithick, goldenrod2002000, mark=*, mark size=\markersize, mark options={solid}]
table {%
6.84008636032105 4
9.30006636032105 4
11.710056360321 4
14.070036360321 4
16.380016360321 4
18.639996360321 4
20.849976360321 4
23.089956360321 4
25.359946360321 4
27.659926360321 4
29.989906360321 4
32.349886360321 4
34.659866360321 4
36.919846360321 4
39.129826360321 4
41.289816360321 4
43.399796360321 4
45.459776360321 4
47.469756360321 4
49.429736360321 4
51.339716360321 4
};
\addplot [semithick, steelblue31119180, mark=*, mark size=\markersize, mark options={solid}]
table {%
7.34422636032104 0.49988
9.80523636032104 1.031227
12.363506360321 1.25123
14.979466360321 1.346801
17.646976360321 1.313773
20.364686360321 1.3104
23.131906360321 1.36288
25.931316360321 1.305675
28.731176360321 1.333971
31.529986360321 1.415505
34.3191263603211 1.661905
37.109316360321 1.896138
39.888286360321 2.238607
42.646336360321 2.721509
45.376556360321 3.342721
48.132416360321 3.837927
50.9279663603211 3.995493
53.727866360321 4.01203
56.527786360321 4.00481
59.327706360321 4.000748
62.127606360321 3.999675
};
\addplot [semithick, goldenrod2002000, mark=*, mark size=\markersize, mark options={solid}]
table {%
6.84008636032105 4
9.30006636032105 4
11.710056360321 4
14.070036360321 4
16.380016360321 4
18.639996360321 4
20.849976360321 4
23.089956360321 4
25.359946360321 4
27.659926360321 4
29.989906360321 4
32.349886360321 4
34.739866360321 4
37.158546360321 4
39.603926360321 4
42.074126360321 4
44.567406360321 4
47.082166360321 4
49.616886360321 4
52.170176360321 4
54.7407363603211 4
};
\addplot [semithick, steelblue31119180, mark=*, mark size=\markersize, mark options={solid}]
table {%
7.34422636032104 0.49988
9.80523636032104 1.031227
12.364356360321 1.241056
14.980206360321 1.339871
17.647466360321 1.389095
20.365176360321 1.392857
23.132046360321 1.324636
25.932026360321 1.313209
28.732026360321 1.314721
31.531546360321 1.366417
34.331546360321 1.36651
37.131006360321 1.311725
39.930076360321 1.383903
42.718966360321 1.633105
45.486116360321 2.060777
48.139566360321 2.954744
50.828466360321 3.735672
53.617056360321 3.9881
56.416886360321 4.017836
59.216846360321 4.007705
62.016806360321 4.00044
};
\addplot [semithick, goldenrod2002000, mark=*, mark size=\markersize, mark options={solid}]
table {%
6.84008636032105 4
9.30006636032105 4
11.783856360321 4
14.289776360321 4
16.816276360321 4
19.361926360321 4
21.925376360321 4
24.438826360321 4
26.902286360321 4
29.315736360321 4
31.679186360321 4
33.992646360321 4
36.256096360321 4
38.469546360321 4
40.633006360321 4
42.746456360321 4
44.809906360321 4
46.823356360321 4
48.786816360321 4
50.700266360321 4
52.563716360321 4
};
\addplot [semithick, steelblue31119180, mark=*, mark size=\markersize, mark options={solid}]
table {%
7.34422636032104 0.49988
9.80523636032104 1.031227
12.364356360321 1.241056
14.980206360321 1.339871
17.647466360321 1.389095
20.365176360321 1.392857
23.132046360321 1.324636
25.932026360321 1.313209
28.732026360321 1.314721
31.531546360321 1.366417
34.331546360321 1.36651
37.131006360321 1.311725
39.930646360321 1.357031
42.728346360321 1.470414
45.523086360321 1.642003
48.314156360321 1.865429
51.099266360321 2.153805
53.878236360321 2.496299
56.636836360321 2.976053
59.374696360321 3.562672
62.146916360321 3.956072
};
\addplot [semithick, goldenrod2002000, mark=*, mark size=\markersize, mark options={solid}]
table {%
6.84008636032105 4
9.30006636032105 4
11.783856360321 4
14.289776360321 4
16.816276360321 4
19.361926360321 4
21.925376360321 4
24.438826360321 4
26.902286360321 4
29.315736360321 4
31.679186360321 4
33.992646360321 4
36.336096360321 4
38.709546360321 4
41.112856360321 4
43.543936360321 4
46.000846360321 4
48.481766360321 4
50.985016360321 4
53.5090463603211 4
56.052396360321 4
};
\addplot [semithick, steelblue31119180, mark=*, mark size=\markersize, mark options={solid}]
table {%
7.34422636032104 0.49988
9.80523636032104 1.031227
12.364356360321 1.241056
14.980206360321 1.339871
17.647466360321 1.389095
20.365176360321 1.392857
23.132046360321 1.324636
25.932006360321 1.310698
28.732006360321 1.310008
31.532006360321 1.314232
34.331966360321 1.330303
37.131466360321 1.383079
39.931376360321 1.405225
42.730936360321 1.355334
45.530586360321 1.311107
48.330576360321 1.316328
51.129746360321 1.384385
53.923546360321 1.570756
56.707986360321 1.865526
59.438616360321 2.484922
62.103016360321 3.345716
};
\addplot [semithick, goldenrod2002000, mark=*, mark size=\markersize, mark options={solid}]
table {%
6.84008636032105 4
9.30006636032105 4
11.783856360321 4
14.289776360321 4
16.816276360321 4
19.361926360321 4
21.925376360321 4
24.505386360321 4
27.100796360321 4
29.710526360321 4
32.333586360321 4
34.969016360321 4
37.554456360321 4
40.089896360321 4
42.575326360321 4
45.010766360321 4
47.396206360321 4
49.731636360321 4
52.017076360321 4
54.252516360321 4
56.4379563603211 4
};
\addplot [semithick, steelblue31119180, mark=*, mark size=\markersize, mark options={solid}]
table {%
7.34422636032104 0.49988
9.80523636032104 1.031227
12.364356360321 1.241056
14.980206360321 1.339871
17.647466360321 1.389095
20.365176360321 1.392857
23.132046360321 1.324636
25.932006360321 1.310698
28.732006360321 1.310008
31.532006360321 1.314232
34.331966360321 1.330303
37.131466360321 1.383079
39.931456360321 1.389185
42.694076360321 1.417529
45.420556360321 1.435437
48.118936360321 1.440853
50.805716360321 1.441977
53.497346360321 1.442092
56.1982063603211 1.442096
58.9061163603211 1.442096
61.620476360321 1.442096
};
\addplot [semithick, goldenrod2002000, mark=*, mark size=\markersize, mark options={solid}]
table {%
6.84008636032105 4
9.30006636032105 4
11.783856360321 4
14.289776360321 4
16.816276360321 4
19.361926360321 4
21.925376360321 4
24.505386360321 4
27.100796360321 4
29.710526360321 4
32.333586360321 4
34.969016360321 4
37.615976360321 4
40.273646360321 4
42.941276360321 4
45.618176360321 4
48.303686360321 4
50.997216360321 4
53.6982063603211 4
56.4061163603211 4
59.120476360321 4
};

\nextgroupplot[
height=\simulationheight,
width=\simulationwidth,
tick align=outside,
tick pos=left,
x grid style={darkgray176},
xmajorgrids,
xmin=42.6271695014572, xmax=124.778732652321,
xtick style={color=black, font=\scriptsize},
y dir=reverse,
y grid style={darkgray176},
ylabel={\(\displaystyle p_y\) [$\SI{}{\meter}$]},
ticklabel style={font=\footnotesize},
ticklabel style={font=\footnotesize},
ymajorgrids,
ymin=-1.5, ymax=5.5,
ytick style={color=black, font=\scriptsize}
]
\node[fill=white, anchor=north east] at (rel axis cs:1,1) {$t = \SI{2.2}{\second}$};
\draw[draw=black,fill=steelblue31119180,opacity=0.5,rotate around={-0.919162462270711:(axis cs:51.1299285888672,1.40506315231323)}] (axis cs:48.6299285888672,0.405063152313232) rectangle (axis cs:53.6299285888672,2.40506315231323);
\draw[draw=black,fill=steelblue31119180,opacity=0.5,rotate around={-0.169920699119637:(axis cs:53.9299011230469,1.41728925704956)}] (axis cs:51.4299011230469,0.417289257049561) rectangle (axis cs:56.4299011230469,2.41728925704956);
\draw[draw=black,fill=steelblue31119180,opacity=0.5,rotate around={1.25386211544714:(axis cs:56.7293853759766,1.3634672164917)}] (axis cs:54.2293853759766,0.363467216491699) rectangle (axis cs:59.2293853759766,2.3634672164917);
\draw[draw=black,fill=steelblue31119180,opacity=0.5,rotate around={0.843015564700295:(axis cs:59.5290489196777,1.32012033462524)}] (axis cs:57.0290489196777,0.320120334625244) rectangle (axis cs:62.0290489196777,2.32012033462524);
\draw[draw=black,fill=steelblue31119180,opacity=0.5,rotate around={0.149041090971385:(axis cs:62.329029083252,1.30920374393463)}] (axis cs:59.829029083252,0.309203743934631) rectangle (axis cs:64.829029083252,2.30920374393463);
\draw[draw=black,fill=steelblue31119180,rotate around={-0.149936179257988:(axis cs:65.129020690918,1.31496560573578)}] (axis cs:62.629020690918,0.314965605735779) rectangle (axis cs:67.629020690918,2.31496560573578);
\draw[draw=black,fill=goldenrod2002000,opacity=0.5] (axis cs:46.3613314628601,3) rectangle (axis cs:51.3613314628601,5);
\draw[draw=black,fill=goldenrod2002000,opacity=0.5] (axis cs:48.9960741996765,3) rectangle (axis cs:53.9960741996765,5);
\draw[draw=black,fill=goldenrod2002000,opacity=0.5] (axis cs:51.5808191299438,3) rectangle (axis cs:56.5808191299438,5);
\draw[draw=black,fill=goldenrod2002000,opacity=0.5] (axis cs:54.1155619621277,3) rectangle (axis cs:59.1155619621277,5);
\draw[draw=black,fill=goldenrod2002000,opacity=0.5] (axis cs:56.6803066730499,2.99999988079071) rectangle (axis cs:61.6803066730499,4.99999988079071);
\draw[draw=black,fill=goldenrod2002000] (axis cs:59.2750477790833,2.99999988079071) rectangle (axis cs:64.2750477790833,4.99999988079071);
\addplot [semithick, steelblue31119180, mark=*, mark size=\markersize, mark options={solid}]
table {%
65.129020690918 1.314966
67.928850690918 1.345619
70.727830690918 1.421146
73.523790690918 1.571601
76.313450690918 1.812039
79.079980690918 2.243658
81.774430690918 3.00518
84.477450690918 3.735746
87.266370690918 3.984545
90.066180690918 4.01645
92.866140690918 4.007333
95.666110690918 4.001288
98.466070690918 3.999848
101.266030690918 3.999835
104.065980690918 3.999952
106.865930690918 3.999997
109.665880690918 4.000003
112.465820690918 4.000001
115.265760690918 4
118.065710690918 4
120.865650690918 4
};
\addplot [semithick, goldenrod2002000, mark=*, mark size=\markersize, mark options={solid}]
table {%
61.775050690918 4
64.319790690918 4
66.814530690918 4
69.259280690918 4
71.654020690918 4
73.998760690918 4
76.293510690918 4
78.538250690918 4
80.732990690918 4
82.877740690918 4
84.972480690918 4
87.017230690918 4
89.011970690918 4
90.956710690918 4
92.851460690918 4
94.696200690918 4
96.490940690918 4
98.235690690918 4
99.930430690918 4
101.575170690918 4
103.169920690918 4
};
\addplot [semithick, steelblue31119180, mark=*, mark size=\markersize, mark options={solid}]
table {%
65.129020690918 1.314966
67.928850690918 1.345619
70.727830690918 1.421146
73.523790690918 1.571601
76.313450690918 1.812039
79.079980690918 2.243658
81.774430690918 3.00518
84.477450690918 3.735746
87.266370690918 3.984545
90.066180690918 4.01645
92.866140690918 4.007333
95.666110690918 4.001288
98.466070690918 3.999848
101.266030690918 3.999835
104.065980690918 3.999952
106.865930690918 3.999997
109.665880690918 4.000003
112.465820690918 4.000001
115.265760690918 4
118.065710690918 4
120.865650690918 4
};
\addplot [semithick, goldenrod2002000, mark=*, mark size=\markersize, mark options={solid}]
table {%
61.775050690918 4
64.319790690918 4
66.814530690918 4
69.259280690918 4
71.654020690918 4
73.998760690918 4
76.293510690918 4
78.538250690918 4
80.732990690918 4
82.877740690918 4
84.972480690918 4
87.017230690918 4
89.091970690918 4
91.196710690918 4
93.331460690918 4
95.496200690918 4
97.690940690918 4
99.915690690918 4
102.170430690918 4
104.455170690918 4
106.769920690918 4
};
\addplot [semithick, steelblue31119180, mark=*, mark size=\markersize, mark options={solid}]
table {%
65.129020690918 1.314966
67.928850690918 1.345619
70.727830690918 1.421146
73.523790690918 1.571601
76.313450690918 1.812039
79.079980690918 2.243658
81.774430690918 3.00518
84.501090690918 3.641841
87.284990690918 3.941723
90.084180690918 4.008815
92.884160690918 4.008168
95.684130690918 4.002244
98.484090690918 4.000127
101.284050690918 3.999856
104.084010690918 3.999936
106.883960690918 3.999989
109.683910690918 4.000001
112.483850690918 4.000001
115.283800690918 4
118.083740690918 4
120.883680690918 4
};
\addplot [semithick, goldenrod2002000, mark=*, mark size=\markersize, mark options={solid}]
table {%
61.775050690918 4
64.319790690918 4
66.814530690918 4
69.259280690918 4
71.654020690918 4
73.998760690918 4
76.293510690918 4
78.618250690918 4
80.972990690918 4
83.357740690918 4
85.771550690918 4
88.212390690918 4
90.603240690918 4
92.944080690918 4
95.234930690918 4
97.475770690918 4
99.666620690918 4
101.807460690918 4
103.898310690918 4
105.939150690918 4
107.930000690918 4
};
\addplot [semithick, steelblue31119180, mark=*, mark size=\markersize, mark options={solid}]
table {%
65.129020690918 1.314966
67.928850690918 1.345619
70.727830690918 1.421146
73.523790690918 1.571601
76.313450690918 1.812039
79.079980690918 2.243658
81.774430690918 3.00518
84.501090690918 3.641841
87.284990690918 3.941723
90.084180690918 4.008815
92.884160690918 4.008168
95.684130690918 4.002244
98.484090690918 4.000127
101.284050690918 3.999856
104.084010690918 3.999936
106.883960690918 3.999989
109.683910690918 4.000001
112.483850690918 4.000001
115.283800690918 4
118.083740690918 4
120.883680690918 4
};
\addplot [semithick, goldenrod2002000, mark=*, mark size=\markersize, mark options={solid}]
table {%
61.775050690918 4
64.319790690918 4
66.814530690918 4
69.259280690918 4
71.654020690918 4
73.998760690918 4
76.293510690918 4
78.618250690918 4
80.972990690918 4
83.357740690918 4
85.771550690918 4
88.212390690918 4
90.678380690918 4
93.167750690918 4
95.678860690918 4
98.210190690918 4
100.760330690918 4
103.327960690918 4
105.911860690918 4
108.510880690918 4
111.123970690918 4
};
\addplot [semithick, steelblue31119180, mark=*, mark size=\markersize, mark options={solid}]
table {%
65.129020690918 1.314966
67.928850690918 1.345619
70.728650690918 1.379033
73.527600690918 1.455893
76.327150690918 1.506118
79.123460690918 1.649861
81.922800690918 1.710363
84.716890690918 1.892177
87.503600690918 2.164691
90.272470690918 2.581031
92.982960690918 3.28337
95.729430690918 3.828233
98.524350690918 3.99664
101.324250690918 4.013234
104.124170690918 4.005133
106.924060690918 4.000756
109.723950690918 3.99984
112.523830690918 3.999877
115.323690690918 3.999968
118.123550690918 3.999999
120.923400690918 4.000003
};
\addplot [semithick, goldenrod2002000, mark=*, mark size=\markersize, mark options={solid}]
table {%
61.775050690918 4
64.319790690918 4
66.882400690918 4
69.461630690918 4
72.056310690918 4
74.665370690918 4
77.287790690918 4
79.860210690918 4
82.382630690918 4
84.855050690918 4
87.277470690918 4
89.649890690918 4
91.972320690918 4
94.244740690918 4
96.467160690918 4
98.639580690918 4
100.762000690918 4
102.834420690918 4
104.856840690918 4
106.829260690918 4
108.751680690918 4
};
\addplot [semithick, steelblue31119180, mark=*, mark size=\markersize, mark options={solid}]
table {%
65.129020690918 1.314966
67.928850690918 1.345619
70.728650690918 1.379033
73.527600690918 1.455893
76.327150690918 1.506118
79.123460690918 1.649861
81.922800690918 1.710363
84.716890690918 1.892177
87.503600690918 2.164691
90.272470690918 2.581031
92.982960690918 3.28337
95.729430690918 3.828233
98.524350690918 3.99664
101.324250690918 4.013234
104.124170690918 4.005133
106.924060690918 4.000756
109.723950690918 3.99984
112.523830690918 3.999877
115.323690690918 3.999968
118.123550690918 3.999999
120.923400690918 4.000003
};
\addplot [semithick, goldenrod2002000, mark=*, mark size=\markersize, mark options={solid}]
table {%
61.775050690918 4
64.319790690918 4
66.882400690918 4
69.461630690918 4
72.056310690918 4
74.665370690918 4
77.287790690918 4
79.860210690918 4
82.382630690918 4
84.855050690918 4
87.277470690918 4
89.649890690918 4
92.052250690918 4
94.482430690918 4
96.938510690918 4
99.418660690918 4
101.921190690918 4
104.444550690918 4
106.987280690918 4
109.548010690918 4
112.125500690918 4
};
\addplot [semithick, steelblue31119180, mark=*, mark size=\markersize, mark options={solid}]
table {%
65.129020690918 1.314966
67.928850690918 1.345619
70.728650690918 1.379033
73.527600690918 1.455893
76.327150690918 1.506118
79.123460690918 1.649861
81.922800690918 1.710363
84.721290690918 1.802279
87.518530690918 1.926536
90.317240690918 2.011497
93.110110690918 2.211263
95.908870690918 2.294566
98.696830690918 2.553911
101.453420690918 3.045065
104.183970690918 3.664826
106.969270690918 3.951287
109.768610690918 4.010373
112.568550690918 4.00793
115.368460690918 4.002029
118.168370690918 4.000066
120.968260690918 3.99978
};
\addplot [semithick, goldenrod2002000, mark=*, mark size=\markersize, mark options={solid}]
table {%
61.775050690918 4
64.319790690918 4
66.882400690918 4
69.461630690918 4
72.056310690918 4
74.665370690918 4
77.287790690918 4
79.922640690918 4
82.569050690918 4
85.226220690918 4
87.893380690918 4
90.569840690918 4
93.196300690918 4
95.772760690918 4
98.299220690918 4
100.775680690918 4
103.202140690918 4
105.578600690918 4
107.905060690918 4
110.181520690918 4
112.407980690918 4
};
\addplot [semithick, steelblue31119180, mark=*, mark size=\markersize, mark options={solid}]
table {%
65.129020690918 1.314966
67.928850690918 1.345619
70.728650690918 1.379033
73.527600690918 1.455893
76.327150690918 1.506118
79.123460690918 1.649861
81.922800690918 1.710363
84.721290690918 1.802279
87.518530690918 1.926536
90.317240690918 2.011497
93.110110690918 2.211263
95.908870690918 2.294566
98.705780690918 2.426138
101.501990690918 2.571806
104.297530690918 2.729792
107.092630690918 2.895227
109.885740690918 3.091587
112.678980690918 3.285987
115.461440690918 3.598941
118.247460690918 3.878352
121.044570690918 4.00495
};
\addplot [semithick, goldenrod2002000, mark=*, mark size=\markersize, mark options={solid}]
table {%
61.775050690918 4
64.319790690918 4
66.882400690918 4
69.461630690918 4
72.056310690918 4
74.665370690918 4
77.287790690918 4
79.922640690918 4
82.569050690918 4
85.226220690918 4
87.893380690918 4
90.569840690918 4
93.254950690918 4
95.948100690918 4
98.648730690918 4
101.356310690918 4
104.070370690918 4
106.790440690918 4
109.516100690918 4
112.246970690918 4
114.982680690918 4
};

\nextgroupplot[
height=\simulationheight,
width=\simulationwidth,
tick align=outside,
tick pos=left,
x grid style={darkgray176},
xmajorgrids,
xmin=52.9080069567757, xmax=135.898600714809,
xtick style={color=black, font=\scriptsize},
y dir=reverse,
y grid style={darkgray176},
ylabel={\(\displaystyle p_y\) [$\SI{}{\meter}$]},
ticklabel style={font=\footnotesize},
ticklabel style={font=\footnotesize},
ymajorgrids,
ymin=-1.5, ymax=5.5,
ytick style={color=black, font=\scriptsize}
]
\node[fill=white, anchor=north east] at (rel axis cs:1,1) {$t = \SI{2.6}{\second}$};
\draw[draw=black,fill=steelblue31119180,opacity=0.5,rotate around={0.149041090971385:(axis cs:62.329029083252,1.30920374393463)}] (axis cs:59.829029083252,0.309203743934631) rectangle (axis cs:64.829029083252,2.30920374393463);
\draw[draw=black,fill=steelblue31119180,opacity=0.5,rotate around={-0.149936179257988:(axis cs:65.129020690918,1.31496560573578)}] (axis cs:62.629020690918,0.314965605735779) rectangle (axis cs:67.629020690918,2.31496560573578);
\draw[draw=black,fill=steelblue31119180,opacity=0.5,rotate around={-0.684529509385759:(axis cs:67.9288558959961,1.34561860561371)}] (axis cs:65.4288558959961,0.345618605613708) rectangle (axis cs:70.4288558959961,2.34561860561371);
\draw[draw=black,fill=steelblue31119180,opacity=0.5,rotate around={-1.77204925114726:(axis cs:70.7276840209961,1.42651402950287)}] (axis cs:68.2276840209961,0.426514029502869) rectangle (axis cs:73.2276840209961,2.42651402950287);
\draw[draw=black,fill=steelblue31119180,opacity=0.5,rotate around={-1.88454887007972:(axis cs:73.5261917114258,1.51800513267517)}] (axis cs:71.0261917114258,0.518005132675171) rectangle (axis cs:76.0261917114258,2.51800513267517);
\draw[draw=black,fill=steelblue31119180,rotate around={-3.69478272066876:(axis cs:76.3209609985352,1.68899190425873)}] (axis cs:73.8209609985352,0.688991904258728) rectangle (axis cs:78.8209609985352,2.68899190425873);
\draw[draw=black,fill=goldenrod2002000,opacity=0.5] (axis cs:56.6803066730499,2.99999988079071) rectangle (axis cs:61.6803066730499,4.99999988079071);
\draw[draw=black,fill=goldenrod2002000,opacity=0.5] (axis cs:59.2750477790833,2.99999988079071) rectangle (axis cs:64.2750477790833,4.99999988079071);
\draw[draw=black,fill=goldenrod2002000,opacity=0.5] (axis cs:61.8197944164276,3.00000011920929) rectangle (axis cs:66.8197944164276,5.00000011920929);
\draw[draw=black,fill=goldenrod2002000,opacity=0.5] (axis cs:64.3145349025726,3.00000011920929) rectangle (axis cs:69.3145349025726,5.00000011920929);
\draw[draw=black,fill=goldenrod2002000,opacity=0.5] (axis cs:66.8392810821533,3) rectangle (axis cs:71.8392810821533,5);
\draw[draw=black,fill=goldenrod2002000] (axis cs:69.3140196800232,2.99999988079071) rectangle (axis cs:74.3140196800232,4.99999988079071);
\addplot [semithick, steelblue31119180, mark=*, mark size=\markersize, mark options={solid}]
table {%
76.3209609985352 1.688992
79.1053609985352 1.984189
81.8748009985352 2.396714
84.5458509985352 3.236639
87.2823909985352 3.829415
90.0771509985352 4.000616
92.8771109985352 4.014557
95.6770909985352 4.005279
98.4770709985352 4.000697
101.277060998535 3.999806
104.077050998535 3.99987
106.877030998535 3.999969
109.677010998535 4
112.476990998535 4.000003
115.276970998535 4.000001
118.076940998535 4
120.876920998535 4
123.676900998535 4
126.476870998535 4
129.276840998535 4
132.076820998535 4
};
\addplot [semithick, goldenrod2002000, mark=*, mark size=\markersize, mark options={solid}]
table {%
71.8140209985352 4
74.2387609985352 4
76.6135009985352 4
78.9382509985352 4
81.2129909985352 4
83.4377409985352 4
85.6124809985352 4
87.7372209985352 4
89.8119709985352 4
91.8367109985352 4
93.8114509985352 4
95.7362009985352 4
97.6109409985352 4
99.4356809985352 4
101.210430998535 4
102.935170998535 4
104.609920998535 4
106.234660998535 4
107.809400998535 4
109.334150998535 4
110.808890998535 4
};
\addplot [semithick, steelblue31119180, mark=*, mark size=\markersize, mark options={solid}]
table {%
76.3209609985352 1.688992
79.1053609985352 1.984189
81.8748009985352 2.396714
84.5458509985352 3.236639
87.2823909985352 3.829415
90.0771509985352 4.000616
92.8771109985352 4.014557
95.6770909985352 4.005279
98.4770709985352 4.000697
101.277060998535 3.999806
104.077050998535 3.99987
106.877030998535 3.999969
109.677010998535 4
112.476990998535 4.000003
115.276970998535 4.000001
118.076940998535 4
120.876920998535 4
123.676900998535 4
126.476870998535 4
129.276840998535 4
132.076820998535 4
};
\addplot [semithick, goldenrod2002000, mark=*, mark size=\markersize, mark options={solid}]
table {%
71.8140209985352 4
74.2387609985352 4
76.6135009985352 4
78.9382509985352 4
81.2129909985352 4
83.4377409985352 4
85.6124809985352 4
87.7372209985352 4
89.8119709985352 4
91.8367109985352 4
93.8114509985352 4
95.7362009985352 4
97.6909409985352 4
99.6756809985352 4
101.690430998535 4
103.735170998535 4
105.809920998535 4
107.914660998535 4
110.049400998535 4
112.214150998535 4
114.408890998535 4
};
\addplot [semithick, steelblue31119180, mark=*, mark size=\markersize, mark options={solid}]
table {%
76.3209609985352 1.688992
79.1053609985352 1.984189
81.8748009985352 2.396714
84.5458509985352 3.236639
87.2823909985352 3.829415
90.0771509985352 4.000616
92.8771109985352 4.014557
95.6770909985352 4.005279
98.4770709985352 4.000697
101.277060998535 3.999806
104.077050998535 3.99987
106.877030998535 3.999969
109.677010998535 4
112.476990998535 4.000003
115.276970998535 4.000001
118.076940998535 4
120.876920998535 4
123.676900998535 4
126.476870998535 4
129.276840998535 4
132.076820998535 4
};
\addplot [semithick, goldenrod2002000, mark=*, mark size=\markersize, mark options={solid}]
table {%
71.8140209985352 4
74.2387609985352 4
76.6135009985352 4
78.9382509985352 4
81.2129909985352 4
83.4377409985352 4
85.6124809985352 4
87.8172209985352 4
90.0519709985352 4
92.3167109985352 4
94.6114509985352 4
96.9362009985352 4
99.2109409985352 4
101.435680998535 4
103.610430998535 4
105.735170998535 4
107.809920998535 4
109.834660998535 4
111.809400998535 4
113.734150998535 4
115.608890998535 4
};
\addplot [semithick, steelblue31119180, mark=*, mark size=\markersize, mark options={solid}]
table {%
76.3209609985352 1.688992
79.1053609985352 1.984189
81.8748009985352 2.396714
84.5458509985352 3.236639
87.2823909985352 3.829415
90.0771509985352 4.000616
92.8771109985352 4.014557
95.6770909985352 4.005279
98.4770709985352 4.000697
101.277060998535 3.999806
104.077050998535 3.99987
106.877030998535 3.999969
109.677010998535 4
112.476990998535 4.000003
115.276970998535 4.000001
118.076940998535 4
120.876920998535 4
123.676900998535 4
126.476870998535 4
129.276840998535 4
132.076820998535 4
};
\addplot [semithick, goldenrod2002000, mark=*, mark size=\markersize, mark options={solid}]
table {%
71.8140209985352 4
74.2387609985352 4
76.6135009985352 4
78.9382509985352 4
81.2129909985352 4
83.4377409985352 4
85.6124809985352 4
87.8172209985352 4
90.0519709985352 4
92.3167109985352 4
94.6114509985352 4
96.9362009985352 4
99.2909409985351 4
101.675680998535 4
104.089500998535 4
106.530340998535 4
108.996330998535 4
111.485690998535 4
113.996800998535 4
116.528140998535 4
119.078280998535 4
};
\addplot [semithick, steelblue31119180, mark=*, mark size=\markersize, mark options={solid}]
table {%
76.3209609985352 1.688992
79.1053609985352 1.984189
81.8978809985352 2.188719
84.6635409985352 2.625846
87.3923309985352 3.253339
90.1370209985351 3.807116
92.9309209985352 3.991406
95.7307409985352 4.013097
98.5305809985352 4.005527
101.330380998535 4.000912
104.130130998535 3.999864
106.929850998535 3.999872
109.729540998535 3.999965
112.529190998535 3.999998
115.328820998535 4.000002
118.128430998535 4.000001
120.928020998535 4
123.727600998535 4
126.527170998535 4
129.326730998535 4
132.126290998535 4
};
\addplot [semithick, goldenrod2002000, mark=*, mark size=\markersize, mark options={solid}]
table {%
71.8140209985352 4
74.2387609985352 4
76.6897709985352 4
79.1652109985352 4
81.6633709985352 4
84.1826609985352 4
86.7216009985352 4
89.2105409985352 4
91.6494809985352 4
94.0384209985352 4
96.3773509985351 4
98.6662909985351 4
100.905230998535 4
103.094170998535 4
105.233110998535 4
107.322050998535 4
109.360990998535 4
111.349920998535 4
113.288860998535 4
115.177800998535 4
117.016740998535 4
};
\addplot [semithick, steelblue31119180, mark=*, mark size=\markersize, mark options={solid}]
table {%
76.3209609985352 1.688992
79.1053609985352 1.984189
81.8978809985352 2.188719
84.6635409985352 2.625846
87.3923309985352 3.253339
90.1370209985351 3.807116
92.9309209985352 3.991406
95.7307409985352 4.013097
98.5305809985352 4.005527
101.330380998535 4.000912
104.130130998535 3.999864
106.929850998535 3.999872
109.729540998535 3.999965
112.529190998535 3.999998
115.328820998535 4.000002
118.128430998535 4.000001
120.928020998535 4
123.727600998535 4
126.527170998535 4
129.326730998535 4
132.126290998535 4
};
\addplot [semithick, goldenrod2002000, mark=*, mark size=\markersize, mark options={solid}]
table {%
71.8140209985352 4
74.2387609985352 4
76.6897709985352 4
79.1652109985352 4
81.6633709985352 4
84.1826609985352 4
86.7216009985352 4
89.2105409985352 4
91.6494809985352 4
94.0384209985352 4
96.3773509985351 4
98.6662909985351 4
100.985230998535 4
103.334170998535 4
105.713110998535 4
108.121520998535 4
110.557350998535 4
113.018660998535 4
115.503690998535 4
118.010760998535 4
120.538340998535 4
};
\addplot [semithick, steelblue31119180, mark=*, mark size=\markersize, mark options={solid}]
table {%
76.3209609985352 1.688992
79.1053609985352 1.984189
81.8978809985352 2.188719
84.6635409985352 2.625846
87.3923309985352 3.253339
90.1370209985351 3.807116
92.9309209985352 3.991406
95.7307409985352 4.013097
98.5305809985352 4.005527
101.330380998535 4.000912
104.130130998535 3.999864
106.929850998535 3.999872
109.729540998535 3.999965
112.529190998535 3.999998
115.328820998535 4.000002
118.128430998535 4.000001
120.928020998535 4
123.727600998535 4
126.527170998535 4
129.326740998535 4
132.126290998535 4
};
\addplot [semithick, goldenrod2002000, mark=*, mark size=\markersize, mark options={solid}]
table {%
71.8140209985352 4
74.2387609985352 4
76.6897709985352 4
79.1652109985352 4
81.6633709985352 4
84.1826609985352 4
86.7216009985352 4
89.2788109985352 4
91.8530209985352 4
94.4430309985352 4
97.0477509985352 4
99.6661309985352 4
102.234510998535 4
104.752890998535 4
107.221280998535 4
109.639660998535 4
112.008040998535 4
114.326430998535 4
116.594810998535 4
118.813190998535 4
120.981570998535 4
};
\addplot [semithick, steelblue31119180, mark=*, mark size=\markersize, mark options={solid}]
table {%
76.3209609985352 1.688992
79.1053609985352 1.984189
81.8978809985352 2.188719
84.6635409985352 2.625846
87.3923309985352 3.253339
90.1370209985351 3.807116
92.9309209985352 3.991406
95.7307409985352 4.013097
98.5305809985352 4.005527
101.330380998535 4.000912
104.130130998535 3.999864
106.929850998535 3.999872
109.729540998535 3.999965
112.529190998535 3.999998
115.328820998535 4.000002
118.128430998535 4.000001
120.928020998535 4
123.727600998535 4
126.527170998535 4
129.326740998535 4
132.126300998535 4
};
\addplot [semithick, goldenrod2002000, mark=*, mark size=\markersize, mark options={solid}]
table {%
71.8140209985352 4
74.2387609985352 4
76.6897709985352 4
79.1652109985352 4
81.6633709985352 4
84.1826609985352 4
86.7216009985352 4
89.2788109985352 4
91.8530209985352 4
94.4430309985352 4
97.0477509985352 4
99.6661309985352 4
102.297230998535 4
104.940140998535 4
107.594060998535 4
110.258200998535 4
112.931850998535 4
115.614350998535 4
118.305070998535 4
121.003440998535 4
123.708920998535 4
};

\nextgroupplot[
height=\simulationheight,
width=\simulationwidth,
tick align=outside,
tick pos=left,
x grid style={darkgray176},
xlabel={\(\displaystyle p_x\) [m]},
xlabel shift={-5pt},
xmajorgrids,
xmin=67.8817291329193, xmax=152.73659691536,
xtick style={color=black, font=\scriptsize},
y dir=reverse,
y grid style={darkgray176},
ylabel={\(\displaystyle p_y\) [$\SI{}{\meter}$]},
ticklabel style={font=\footnotesize},
ticklabel style={font=\footnotesize},
ymajorgrids,
ymin=-1.5, ymax=5.5,
ytick style={color=black, font=\scriptsize}
]
\node[fill=white, anchor=north east] at (rel axis cs:1,1) {$t = \SI{3.2}{\second}$};
\draw[draw=black,fill=steelblue31119180,opacity=0.5,rotate around={-6.33390420982504:(axis cs:79.1053619384766,1.98418939113617)}] (axis cs:76.6053619384766,0.984189391136169) rectangle (axis cs:81.6053619384766,2.98418939113617);
\draw[draw=black,fill=steelblue31119180,opacity=0.5,rotate around={-9.95075660649807:(axis cs:81.8664321899414,2.44945287704468)}] (axis cs:79.3664321899414,1.44945287704468) rectangle (axis cs:84.3664321899414,3.44945287704468);
\draw[draw=black,fill=steelblue31119180,opacity=0.5,rotate around={-17.9170649288988:(axis cs:84.5429153442383,3.27192163467407)}] (axis cs:82.0429153442383,2.27192163467407) rectangle (axis cs:87.0429153442383,4.27192163467407);
\draw[draw=black,fill=steelblue31119180,opacity=0.5,rotate around={-10.9491910980496:(axis cs:87.2849044799805,3.83893442153931)}] (axis cs:84.7849044799805,2.83893442153931) rectangle (axis cs:89.7849044799805,4.83893442153931);
\draw[draw=black,fill=steelblue31119180,opacity=0.5,rotate around={-2.43501985538726:(axis cs:90.080192565918,4.0013222694397)}] (axis cs:87.580192565918,3.0013222694397) rectangle (axis cs:92.580192565918,5.0013222694397);
\draw[draw=black,fill=steelblue31119180,rotate around={0.000409750537588144:(axis cs:92.8801574707031,4.01402854919434)}] (axis cs:90.3801574707031,3.01402854919434) rectangle (axis cs:95.3801574707031,5.01402854919434);
\draw[draw=black,fill=goldenrod2002000,opacity=0.5] (axis cs:71.7387685775757,3.00000011920929) rectangle (axis cs:76.7387685775757,5.00000011920929);
\draw[draw=black,fill=goldenrod2002000,opacity=0.5] (axis cs:74.1135082244873,3) rectangle (axis cs:79.1135082244873,5);
\draw[draw=black,fill=goldenrod2002000,opacity=0.5] (axis cs:76.4382548332214,2.99999994039536) rectangle (axis cs:81.4382548332214,4.99999994039536);
\draw[draw=black,fill=goldenrod2002000,opacity=0.5] (axis cs:78.7129998207092,2.99999991059303) rectangle (axis cs:83.7129998207092,4.99999991059303);
\draw[draw=black,fill=goldenrod2002000,opacity=0.5] (axis cs:80.937744140625,3.00000001105946) rectangle (axis cs:85.937744140625,5.00000001105946);
\draw[draw=black,fill=goldenrod2002000] (axis cs:83.1124868392944,3.000000230968) rectangle (axis cs:88.1124868392944,5.000000230968);
\addplot [semithick, steelblue31119180, mark=*, mark size=\markersize, mark options={solid}]
table {%
92.8801574707031 4.014029
95.6801374707031 4.005026
98.4801174707031 4.000649
101.280097470703 3.99981
104.080057470703 3.999875
106.880017470703 3.99997
109.679977470703 4
112.479917470703 4.000002
115.279857470703 4.000001
118.079787470703 4
120.879717470703 4
123.679637470703 4
126.479557470703 4
129.279467470703 4
132.079377470703 4
134.879287470703 4
137.679187470703 4
140.479087470703 4
143.278987470703 4
146.078887470703 4
148.878787470703 4
};
\addplot [semithick, goldenrod2002000, mark=*, mark size=\markersize, mark options={solid}]
table {%
85.6124874707031 4
87.7372274707031 4
89.8119774707031 4
91.8367174707031 4
93.8114574707031 4
95.7362074707031 4
97.6109474707031 4
99.4356874707031 4
101.210437470703 4
102.935177470703 4
104.609917470703 4
106.234667470703 4
107.809407470703 4
109.334157470703 4
110.808897470703 4
112.233637470703 4
113.608387470703 4
114.933127470703 4
116.207867470703 4
117.432617470703 4
118.607357470703 4
};
\addplot [semithick, steelblue31119180, mark=*, mark size=\markersize, mark options={solid}]
table {%
92.8801574707031 4.014029
95.6801374707031 4.005026
98.4801174707031 4.000649
101.280097470703 3.99981
104.080057470703 3.999875
106.880017470703 3.99997
109.679977470703 4
112.479917470703 4.000002
115.279857470703 4.000001
118.079787470703 4
120.879717470703 4
123.679637470703 4
126.479557470703 4
129.279467470703 4
132.079377470703 4
134.879287470703 4
137.679187470703 4
140.479087470703 4
143.278987470703 4
146.078887470703 4
148.878787470703 4
};
\addplot [semithick, goldenrod2002000, mark=*, mark size=\markersize, mark options={solid}]
table {%
85.6124874707031 4
87.7372274707031 4
89.8119774707031 4
91.8367174707031 4
93.8114574707031 4
95.7362074707031 4
97.6109474707031 4
99.4356874707031 4
101.210437470703 4
102.935177470703 4
104.609917470703 4
106.234667470703 4
107.889407470703 4
109.574157470703 4
111.288897470703 4
113.033637470703 4
114.808387470703 4
116.613127470703 4
118.447867470703 4
120.312617470703 4
122.207357470703 4
};
\addplot [semithick, steelblue31119180, mark=*, mark size=\markersize, mark options={solid}]
table {%
92.8801574707031 4.014029
95.6801374707031 4.005026
98.4801174707031 4.000649
101.280097470703 3.99981
104.080057470703 3.999875
106.880017470703 3.99997
109.679977470703 4
112.479917470703 4.000002
115.279857470703 4.000001
118.079787470703 4
120.879717470703 4
123.679637470703 4
126.479557470703 4
129.279467470703 4
132.079377470703 4
134.879287470703 4
137.679187470703 4
140.479087470703 4
143.278987470703 4
146.078887470703 4
148.878787470703 4
};
\addplot [semithick, goldenrod2002000, mark=*, mark size=\markersize, mark options={solid}]
table {%
85.6124874707031 4
87.7372274707031 4
89.8119774707031 4
91.8367174707031 4
93.8114574707031 4
95.7362074707031 4
97.6109474707031 4
99.5156874707031 4
101.450437470703 4
103.415177470703 4
105.409917470703 4
107.434667470703 4
109.409407470703 4
111.334157470703 4
113.208897470703 4
115.033637470703 4
116.808387470703 4
118.533127470703 4
120.207867470703 4
121.832617470703 4
123.407357470703 4
};
\addplot [semithick, steelblue31119180, mark=*, mark size=\markersize, mark options={solid}]
table {%
92.8801574707031 4.014029
95.6801374707031 4.005026
98.4801174707031 4.000649
101.280097470703 3.99981
104.080057470703 3.999875
106.880017470703 3.99997
109.679977470703 4
112.479917470703 4.000002
115.279857470703 4.000001
118.079787470703 4
120.879717470703 4
123.679637470703 4
126.479557470703 4
129.279467470703 4
132.079377470703 4
134.879287470703 4
137.679187470703 4
140.479087470703 4
143.278987470703 4
146.078887470703 4
148.878787470703 4
};
\addplot [semithick, goldenrod2002000, mark=*, mark size=\markersize, mark options={solid}]
table {%
85.6124874707031 4
87.7372274707031 4
89.8119774707031 4
91.8367174707031 4
93.8114574707031 4
95.7362074707031 4
97.6109474707031 4
99.5156874707031 4
101.450437470703 4
103.415177470703 4
105.409917470703 4
107.434667470703 4
109.489407470703 4
111.574157470703 4
113.688897470703 4
115.833637470703 4
118.008387470703 4
120.213127470703 4
122.447867470703 4
124.712617470703 4
127.007357470703 4
};
\addplot [semithick, steelblue31119180, mark=*, mark size=\markersize, mark options={solid}]
table {%
92.8801574707031 4.014029
95.6801374707031 4.005026
98.4801174707031 4.000649
101.280097470703 3.99981
104.080077470703 3.999875
106.880057470703 3.99997
109.680037470703 4
112.480017470703 4.000002
115.279987470703 4.000001
118.079957470703 4
120.879927470703 4
123.679897470703 4
126.479867470703 4
129.279827470703 4
132.079797470703 4
134.879757470703 4
137.679717470703 4
140.479677470703 4
143.279637470703 4
146.079597470703 4
148.879557470703 4
};
\addplot [semithick, goldenrod2002000, mark=*, mark size=\markersize, mark options={solid}]
table {%
85.6124874707031 4
87.7372274707031 4
89.8919774707031 4
92.0767174707031 4
94.2914574707031 4
96.5362074707031 4
98.8109474707031 4
101.035687470703 4
103.210437470703 4
105.335177470703 4
107.409917470703 4
109.434667470703 4
111.409407470703 4
113.334157470703 4
115.208897470703 4
117.033637470703 4
118.808387470703 4
120.533127470703 4
122.207867470703 4
123.832617470703 4
125.407357470703 4
};
\addplot [semithick, steelblue31119180, mark=*, mark size=\markersize, mark options={solid}]
table {%
92.8801574707031 4.014029
95.6801374707031 4.005026
98.4801174707031 4.000649
101.280097470703 3.99981
104.080077470703 3.999875
106.880057470703 3.99997
109.680037470703 4
112.480017470703 4.000002
115.279987470703 4.000001
118.079957470703 4
120.879927470703 4
123.679897470703 4
126.479867470703 4
129.279827470703 4
132.079797470703 4
134.879757470703 4
137.679717470703 4
140.479677470703 4
143.279637470703 4
146.079597470703 4
148.879557470703 4
};
\addplot [semithick, goldenrod2002000, mark=*, mark size=\markersize, mark options={solid}]
table {%
85.6124874707031 4
87.7372274707031 4
89.8919774707031 4
92.0767174707031 4
94.2914574707031 4
96.5362074707031 4
98.8109474707031 4
101.035687470703 4
103.210437470703 4
105.335177470703 4
107.409917470703 4
109.434667470703 4
111.489407470703 4
113.574157470703 4
115.688897470703 4
117.833637470703 4
120.008387470703 4
122.213127470703 4
124.447867470703 4
126.712617470703 4
129.007357470703 4
};
\addplot [semithick, steelblue31119180, mark=*, mark size=\markersize, mark options={solid}]
table {%
92.8801574707031 4.014029
95.6801374707031 4.005026
98.4801174707031 4.000649
101.280097470703 3.99981
104.080077470703 3.999875
106.880057470703 3.99997
109.680037470703 4
112.480017470703 4.000002
115.279987470703 4.000001
118.079957470703 4
120.879927470703 4
123.679897470703 4
126.479867470703 4
129.279827470703 4
132.079797470703 4
134.879757470703 4
137.679717470703 4
140.479677470703 4
143.279637470703 4
146.079597470703 4
148.879557470703 4
};
\addplot [semithick, goldenrod2002000, mark=*, mark size=\markersize, mark options={solid}]
table {%
85.6124874707031 4
87.7372274707031 4
89.8919774707031 4
92.0767174707031 4
94.2914574707031 4
96.5362074707031 4
98.8109474707031 4
101.115687470703 4
103.450437470703 4
105.815177470703 4
108.209917470703 4
110.633037470703 4
113.006147470703 4
115.329257470703 4
117.602367470703 4
119.825477470703 4
121.998587470703 4
124.121697470703 4
126.194817470703 4
128.217927470703 4
130.191037470703 4
};
\addplot [semithick, steelblue31119180, mark=*, mark size=\markersize, mark options={solid}]
table {%
92.8801574707031 4.014029
95.6801374707031 4.005026
98.4801174707031 4.000649
101.280097470703 3.99981
104.080077470703 3.999875
106.880057470703 3.99997
109.680037470703 4
112.480017470703 4.000002
115.279987470703 4.000001
118.079957470703 4
120.879927470703 4
123.679897470703 4
126.479867470703 4
129.279827470703 4
132.079797470703 4
134.879757470703 4
137.679717470703 4
140.479677470703 4
143.279637470703 4
146.079597470703 4
148.879557470703 4
};
\addplot [semithick, goldenrod2002000, mark=*, mark size=\markersize, mark options={solid}]
table {%
85.6124874707031 4
87.7372274707031 4
89.8919774707031 4
92.0767174707031 4
94.2914574707031 4
96.5362074707031 4
98.8109474707031 4
101.115687470703 4
103.450437470703 4
105.815177470703 4
108.209917470703 4
110.633037470703 4
113.082527470703 4
115.556557470703 4
118.053407470703 4
120.571467470703 4
123.109277470703 4
125.665427470703 4
128.238657470703 4
130.827757470703 4
133.431617470703 4
};
\end{groupplot}

\end{tikzpicture}

%% file: content/conclusion_and_future_work.tex
We presented a control framework based on \ac{SMPC} for lane changing problem in autonomous driving, where the interaction was accounted for by modeling the behavior of surrounding vehicles using a stochastic policy with state-dependent distribution. 
We present an approximate reformulation of the resulting nonlinear, chance-constrained optimal control problem and presented a simple moving horizon scheme
for learning the parameters describing the driving style of an adjacent driver.
Qualitative and quantitative analysis shows that learning the distribution helps to reduce conservatism of the resulting MPC controller.
